\documentclass{gtart}

\def\ifplaintex{\expandafter\ifx\csname documentclass\endcsname\relax}

\def\gtp{{\mathsurround=0pt\it $\cal G\mskip-2mu$eometry \&\ 
$\cal T\!\!$opology $\cal P\!$ublications}}  

\def\recd{{\small Received:\qua\receiveddate\ifx\reviseddate\relax
\else\qquad Revised:\qua\reviseddate\fi\par}} 

\def\lognumber#1{\def\thelognumber{#1}}
\def\volumenumber#1{\def\thevolumenumber{#1}}
\def\volumeyear#1{\def\thevolumeyear{#1}}
\def\papernumber#1{\def\thepapernumber{#1}}
\def\pagenumbers#1#2{\def\startpage{#1}\def\finishpage{#2}}
\def\published#1{\def\publishdate{#1}}

\def\received#1{\def\receiveddate{#1}}
\def\revised#1{\def\reviseddate{#1}}
\def\accepted#1{\def\accepteddate{#1}}
\def\asciititle#1{\def\theasciititle{#1}}

\def\asciiaddress#1{\def\theasciiaddress{#1}}

\long\def\asciiabstract#1{\long\def\theasciiabstract{#1}}
\def\asciikeywords#1{\def\theasciikeywords{#1}}

\let\\\par\let\thelognumber\relax\let\thevolumenumber\relax
\let\thepapernumber\relax\let\thevolumeyear\relax\let\startpage\relax
\let\finishpage\relax\let\publishdate\relax\let\receiveddate\relax
\let\reviseddate\relax\let\accepteddate\relax\let\theasciititle\relax
\let\theasciiauthors\relax\let\theasciiaddress\relax
\let\theasciiabstract\relax\let\theasciikeywords\relax

\let\theasciiemail\relax

\ifplaintex
\font\logobig=cmssbx10 scaled 3836
\font\logomed=cmssbx10 scaled 2557
\else
\font\logobig=cmssbx10 scaled 4200
\font\logomed=cmssbx10 scaled 2800
\fi

\long\def\makeagttitle{   
\count0=\startpage
\agt\hfill      
\hbox to 45truept{\vbox to 0pt{\vglue -13truept{\logomed A\kern -.37em{\logobig 
T}\kern -.38em G}\vss}\hss}
\break
{\small Volume \thevolumenumber\ (\thevolumeyear)
\startpage--\finishpage\nl
Published: \publishdate}

\vglue .25truein

{\parskip=0pt\leftskip 0pt plus
1fil\def\\{\par\smallskip}{\Large\bf\thetitle}\par\medskip} \vglue
0.05truein

{\parskip=0pt\leftskip 0pt plus 1fil\def\\{\par}{\sc\theauthors}
\par\medskip}%
 
\vglue 0.03truein 

{\small\leftskip 25truept\rightskip 25truept{\bf Abstract}\stdspace\theabstract

{\bf AMS Classification}\stdspace\theprimaryclass
\ifx\thesecondaryclass\relax\else; \thesecondaryclass\fi\par
{\bf Keywords}\stdspace \thekeywords\par}\vglue 7truept

}   

\ifplaintex
\font\phead=cmsl9 scaled 950
\font\pnum=cmbx10 scaled 913
\font\pfoot=cmsl9 scaled 950
\headline{\vbox to 0pt{\vskip -4.5mm\line{\small\phead\ifnum
\count0=\startpage ISSN 1472-2739 (on-line) 1472-2747 (printed)
\hfill {\pnum\folio}\else\ifodd\count0\def\\{ }%
\ifx\theshorttitle\relax\thetitle\else\theshorttitle\fi\hfill{\pnum\folio}
\else\def\\{ and }{\pnum\folio}\hfill\ifx\theshortauthors\relax\theauthors
\else\theshortauthors\fi\fi\fi}\vss}}
\footline{\vbox to 0pt{\vglue 0mm\line{\small\pfoot\ifnum\count0=\startpage
\copyright\ \gtp\hfill\else
\agt, Volume \thevolumenumber\ (\thevolumeyear)\hfill\fi}\vss}}
\else
\font=cmsl9 scaled 1050
\font\lnum=cmbx10 
\font=cmsl9 scaled 1050
\makeatletter
\def\@oddhead{{\small\lhead\ifnum\count0=\startpage ISSN 1472-2739 
(on-line) 1472-2747 (printed)\hfill {\lnum\number\count0}\else\ifodd\count0
\def\\{ }\ifx\theshorttitle\relax \thetitle \else\theshorttitle\fi\hfill
{\lnum\number\count0}\else\def\\{ and }{\lnum\number\count0}
\hfill\ifx\theshortauthors\relax 
\theauthors\else\theshortauthors\fi\fi\fi}}\def\@evenhead{\@oddhead}
\def\@oddfoot{\small\lfoot\ifnum\count0=\startpage\copyright\ \gtp\hfill\else
\agt, Volume \thevolumenumber\ (\thevolumeyear)\hfill\fi}
\def\@evenfoot{\@oddfoot}
\makeatother
\fi
\let\maketitlepage\makeagttitle
\let\makeshorttitle\maketitlepage
\let\maketitle\maketitlepage

\newwrite\gtoutfile
\long\gdef\makeheadfile{  
{\def\\{, }\def\s{ }
\immediate\openout\gtoutfile head.xxx
\immediate\write\gtoutfile{To: math@arxiv.org}
\immediate\write\gtoutfile{Subject: put OR rep NNNNN:ppppp}
\immediate\write\gtoutfile{--text follows this line--}
\immediate\write\gtoutfile{Proxy-for: \ifx\theasciiauthors\relax
\theauthors\else\theasciiauthors\fi\s<\ifx\theasciiemail\relax\theemail\else\theasciiemail\fi>}
\immediate\write\gtoutfile{\noexpand\\}
\immediate\write\gtoutfile{Authors: \ifx\theasciiauthors\relax
\theauthors\else\theasciiauthors\fi}
{\def\\{ }\immediate\write\gtoutfile{Title: \ifx\theasciititle\relax
\thetitle\else\theasciititle\fi}}
\immediate\write\gtoutfile{Subj-class: GT or SG, GR etc}
\immediate\write\gtoutfile{MSC-class: \theprimaryclass\ifx\thesecondaryclass\relax\else, \thesecondaryclass\fi}
\immediate\write\gtoutfile{Journal-ref: Algebr. Geom. Topol. \thevolumenumber\s
(\thevolumeyear) \startpage-\finishpage}
\immediate\write\gtoutfile{Comments: Published by Algebraic and
Geometric Topology at}
\immediate\write\gtoutfile{\s\s\s  http://www.maths.warwick.ac.uk/agt/AGTVol\thevolumenumber/agt-\thevolumenumber-\thepapernumber.abs.html}
\immediate\write\gtoutfile{\noexpand\\}
\immediate\write\gtoutfile{}
\ifx\theasciiabstract\relax
\immediate\write\gtoutfile{\theabstract}\else
\immediate\write\gtoutfile{\theasciiabstract}\fi
\immediate\write\gtoutfile{}
\immediate\write\gtoutfile{\noexpand\\}
\immediate\write\gtoutfile{}
\immediate\closeout\gtoutfile}}  

\def\maketitlepage{\makeagttitle\makeheadfile}
\let\makeshorttitle\maketitlepage
\let\maketitle\maketitlepage

\lognumber{2}
\volumenumber{3}
\volumeyear{2003}
\papernumber{2}
\published{27 January 2002}
\pagenumbers{33}{87}
\received{22 July 2002}
\revised{27 November 2002}
\accepted{10 January 2003}

\usepackage{amsmath,amssymb,graphicx}

\newcommand{\be}{\begin{equation}}
\newcommand{\ee}{\end{equation}}
\newcommand{\ba}{\begin{array}}
\newcommand{\ea}{\end{array}}
\newcommand{\bea}{\begin{eqnarray}}
\newcommand{\eea}{\end{eqnarray}}
\newcommand{\beq}{\begin{eqnarray*}}
\newcommand{\eeq}{\end{eqnarray*}}

\newcommand{\qbi}[2]{\left[\begin{array}{c}#1\\ #2\end{array}\right]}

\newcommand{\fig}[1]{\includegraphics[scale=1]{#1}}

\newcommand{\1}{\mbox {${\bf 1}$}}
\newcommand{\Z}{\mbox {${\cal Z}$}}
\newcommand{\h}{\mbox {${\cal H}$}}

\newcommand{\T}{\mbox{$\cal T$}}
\newcommand{\tr}{\mbox{tr}}

\swapnumbers
\newtheorem{pppar}{}[section]
\swapnumbers
\newtheorem{lemma}[pppar]{Lemma}
\newtheorem{theo}[pppar]{Theorem}
\newtheorem{prop}[pppar]{Proposition}
\newtheorem{corr}[pppar]{Corollary}
\newtheorem{conj}[pppar]{Conjecture}
\theoremstyle{definition}
\newtheorem{dfn}[pppar]{Definition}

\newenvironment{ssec}{\begin{pppar} \rm}{\end{pppar}}

\begin{document}

    \title[HKR--type invariants of 4--thickenings]{HKR--type invariants 
    of 4--thickenings\\of 2--dimensional CW complexes}

\asciititle{HKR-type invariants of 4-thickenings of 2-dimensional CW complexes}

\authors{Ivelina Bobtcheva\\Maria Grazia Messia}
\address{Dipartimento di Scienze Matematiche, Universit\`a di 
Ancona\\Via Brece Bianche 1, 60131, Ancona, Italy}
\asciiaddress{Dipartimento di Scienze Matematiche, Universita di 
Ancona\\Via Brece Bianche 1, 60131, Ancona, Italy}

\email{bobtchev@dipmat.unian.it}

\begin{abstract}
The HKR (Hennings--Kauffman--Radford) framework is used to construct
invariants of 4--thickenings of 2--dimensional CW complexes under
2--deformations (1-- and 2-- handle slides and creations and
cancellations of 1--2 handle pairs).  The input of the invariant is a
finite dimensional unimodular ribbon Hopf algebra $A$ and an element
in a quotient of its center, which determines a trace function on
$A$. We study the subset $\T^4$ of trace elements which define
invariants of 4--thickenings under 2--deformations.  In $\T^4$ two
subsets are identified : $\T^3\subset \T^4$, which produces invariants
of 4--thickenings normalizable to invariants of the boundary, and
$\T^2\subset \T^4$, which produces invariants of 4--thickenings
depending only on the 2--dimensional spine and the second Whitney
number of the 4--thickening.  The case of the quantum $sl(2)$ is
studied in details.  We conjecture that $sl(2)$ leads to four
HKR--type invariants and describe the corresponding trace elements.
Moreover, the fusion algebra of the semisimple quotient of the
category of representations of the quantum $sl(2)$ is identified as a
subalgebra of a quotient of its center.
\end{abstract}

\asciiabstract{
The HKR (Hennings-Kauffman-Radford) framework is used to construct
invariants of 4-thickenings of 2--imensional CW complexes under
2-deformations (1- and 2- handle slides and creations and
cancellations of 1-2 handle pairs).  The input of the invariant is a
finite dimensional unimodular ribbon Hopf algebra A and an element in
a quotient of its center, which determines a trace function on A.  We
study the subset T^4 of trace elements which define invariants of
4-thickenings under 2-deformations.  In ^4 two subsets are identified
: T^3, which produces invariants of 4-thickenings normalizable to
invariants of the boundary, and T^2, which produces invariants of
4-thickenings depending only on the 2-dimensional spine and the second
Whitney number of the 4-thickening.  The case of the quantum sl(2) is
studied in details.  We conjecture that sl(2) leads to four HKR-type
invariants and describe the corresponding trace elements.  Moreover,
the fusion algebra of the semisimple quotient of the category of
representations of the quantum sl(2) is identified as a subalgebra of
a quotient of its center.}

\asciikeywords{Hennings' invariant, Hopf algebras, CW complexes, 4-thickenings}
\maketitle


\section{Introduction}

\begin{ssec} \label{acmoves}
   The (generalized) Andrews--Curtis conjecture \cite{AC}
asserts that any simple homotopy equivalence of 2--complexes can be
obtained by deformation through
2--complexes (expansions
and collapses of disks of dimension at most two and changing the 
attaching maps of the 2--cells by homotopy), to which we refer
here as a
2--deformation. This conjecture is expected to be false and different
proposals for counterexamples have been made, but there seem to be
a lack of tools for actually detecting them as such. An extensive
reference for all the problems connected with the Andrews--Curtis conjecture
is \cite{SAM}.

To any 2--dimensional CW complex $P$, there corresponds a presentation of
its fundamental group, which can be obtained by selecting a vertex as a
base point $b$ and a
spanning tree $T$ in the one-skeleton $P_1$ on the complex. Then any 
1--cell $x_i$ which is not in $T$, with a choice of orientation
 determines an element in $\pi _1(P_1,b)$
and the attaching map of any 2--cell defines a word $R_j$ in the $x_i$'s
which represents a trivial element in $\pi _1(P,b)$.
The  presentation of $\pi _1(P,b)$ obtained in this way,
$\hat{P}=\langle x_1,x_2,\dots ,x_n\mid R_1,R_2,\dots ,R_m\rangle $, depends on the choices
made, but this dependence can be explicitly described. In  \cite{SAM}
(theorem 2.4), it is shown that
the correspondence $P\rightarrow \hat{P}$ induces a bijection between the
2--deformation types of connected 2--dimensional
CW complexes  and the equivalence classes of
finite presentations under the following moves:
\begin{itemize}
\item[(i)] The places of $R_1$ and $R_s$ are interchanged;
\item[(ii)] $R_1$ is replaced with $gR_1g^{-1}$, where $g$ is any
             element in the group, or the reverse of such a move;
\item[(iii)] $R_1$ is replaced with $R_1^{-1}$;
\item[(iv)] $R_1$ is replaced with $R_1R_2$;
\item[(v)] Adding of an additional generator $y$ and an additional
             relator $y R$, where $R$ is any word in the $x_i$'s, or
             the reverse of such a move;
\end{itemize}
We will refer to these six operations as AC--moves and, hopefully
without causing confusion, changing a presentation with a sequence of
AC--moves will be called again a 2--deformation of this presentation.
The inverse $\hat{P}\rightarrow P$ of the bijection above is obtained by
taking one-point union of $n$ circles and attaching on them $m$
2--cells as described by the relations.

If two complexes $X$ and $Y$ are simple
homotopy equivalent, then for some $k$ there exists a 2--deformation
from the one-point union of
$X$ with $k$ copies of $S^{2}$ to  the one-point union of
$Y$ with $k$ copies of $S^{2}$. In particular, if an invariant of
2--complexes under a 2--deformation is
multiplicative under one point union, in order to have some hope of
detecting a counterexample
of the AC--conjecture, its value on $S^2$ should not
be a unit. Since, using the correspondence above, we will talk
instead
about invariants of presentations under the AC--moves, a multiplicative
invariant would be considered potentially interesting for the AC--conjecture if
its value for $\langle \emptyset \mid 1\rangle $ is not a unit.

Such invariants were introduced by
Quinn in
\cite{Q:lectures} and studied in
\cite{B}. The input for their construction is a finite
semisimple {\it symmetric\/} monoidal
category, which is taken to be one of the Lie families
described by Gelfand and Kazhdan in 
\cite{GK}, obtained as subquotients of mod $p$
representations of
simple Lie algebras. Unfortunately, extensive numerical study of
Quinn's invariants (described in \cite{wp})
indicated that, in all numerically generated
examples, the invariants come from
a  representation of the free group on the generators into a
subgroup of
$GL_{N}(Z/p)$ for some $N$,  and in this representation every word has order
$p$. Consequently,  it was shown in  \cite{Klaus} that any invariant
possessing this property can't detect counterexamples to the
AC--conjecture.
  \end{ssec}

In the present work we use the framework of Hennings--Kauffman--Radford (HKR)
\cite{Hennings,KR} to construct
invariants of 4--dimensional thickenings of 2--compl\-exes under
2--deformations, i.e.\ 1-- and 2-- handle  slides and creations or
cancellations of 1--2 handle pairs.
The construction is based on a presentation of a 4--thickening by a
framed link in $S^3$ (where the 1--handles are described by dotted
components) and the input data is a finite
dimensional
unimodular ribbon Hopf algebra  and an element in a quotient of its
center which determines  a trace function on the algebra.

As Hennings points out, any trace function on the algebra, and
therefore any trace element,
leads to an invariant of links, but very few trace elements lead to
invariants of links which are also invariants under the band-connected
sum of two link components (corresponding to 2--handle slides).
Let $\T_s$ be the subset of these special trace elements. Then
$\T_s$ contains always at least two elements which are 1 (one) and the
algebra integral $\Lambda$. Moreover, when the Hopf algebra is the finite
dimensional quantum enveloping algebra at root of unity of some simple
Lie group,  $\T_s$ contains at least one more element $z_{RT}$ which
corresponds to the Reshetikhin--Turaev invariant. This fact was first
observed by Hennings, and then, for the quantum $sl(2)$,  $z_{RT}$ was
made explicit by
Kerler in \cite{K:mcg} (for completeness, in the appendix we present
the derivation of $z_{RT}$). 
In an analogous way (though it won't be done here), one can see
that Quinn's invariant can be derived in the HKR--framework from a
triangular Hopf algebra over $Z/p$ and a central element $z_Q\neq 1$
in it. Moreover, the invariant corresponding to 1 is less
interesting than Quinn's invariant.
These facts imply that it is
important not to restrict to the trace function corresponding to 1
 (as it is done in \cite{Hennings,KR}), and rise the question 
what is the possible
relationship between the different invariants derived from
the same Hopf algebra.  To answer this question,  one needs to  study  the
structure of $\T_s$ and we hope that the present work sets the framework for 
such study.

In particular, we determine  a subset $\T$ of trace elements which 
 lead to  invariants of 4--thickenings under 2--handle slides, i.e.\ $\T\subset \T_s$.
By adding the requirement for invariance
under 1--2 handle cancellations,
inside $\T$, it is described a subset  $\T^4$
of trace elements which lead to
invariants of 4--thickenings under 2--deformations.
Then we study when the invariant of a 4--thickening reduces to un 
invariant of   the boundary and when it reduces to an invariant of
the spine. This leads to the description in $\T^4$
of two  subsets:
\begin{itemize}
    \item $\T^3\subset \T^4$, whose elements lead to invariants which
    factor as a product of a 3-manifold invariant and a 
    multiplicative invariant which
depends on  the signature and the Euler characteristic of
the 4--thickening, and 
   \item $\T^2\subset \T^4$ whose elements  lead to invariants which 
   only depend on the 2--dimensional
spine and the second Whitney number of the 4--thickening.
\end{itemize}

The definition of $\T$ allows to make some interesting
conclusions about its structure. In particular, $\T$ carries two
different monoidal structures and it is invariant under the action of
the ${\cal S}$ operator
  defined in (2.55) of \cite{K:mcg}.
But for now we don't know a practical way of calculating
the elements of $\T$ for a given algebra
and this is quite unsatisfactory. The only partial
remedy we can offer, is that by weakening ``slightly'' the defining conditions
on $\T$, one can define a subset $\T_Z$, containing $\T$,
such that its elements
 are relatively easy to determine since the calculations are entirely
restricted to the center of the algebra. We make this calculation
explicit for the case of the quantum $sl(2)$ and show that in this
case $\T_Z$
consists of 4 elements, three of which are exactly $1,\,\Lambda$ and
$z_{RT}$. This fact leads to the conjecture that $\T_Z=\T$ for the quantum
$sl(2)$. Under this assumption we show that the invariant
corresponding to the forth element
in $\T_Z$ is  the ratio of the Hennings and the
Reshetikhin--Turaev invariants.

The paper is organized as follows. In section 2 we present the main 
definitions and results.
Section 3 contains some notations and preliminaries on Hopf algebras.
Section 4 is dedicated to the study  of the structure of $\T$.
Section 5 introduces the notion of K--links and K--tangles.
Section 6  defines the invariant
of 4--thickenings and shows that, when the trace element is in $\T^2$,
the invariant depends only on the two
dimensional spine  of the 4--thickening and its second Whitney number.
Section 7 studies the reducibility of the
invariant to a 3--manifold invariant
and section 8 illustrates the
construction with two examples: the case of a group algebra and the
case of the quantum  $sl(2)$. At the end we list some open
questions. In the appendix, always for the quantum $sl(2)$, we show
that the Reshetikhin--Turaev invariant is a HKR--type invariant and
calculate the corresponding trace element.

 {\bf Acknowledgements}\qua We want to thank   Thomas 
 Kerler,  Frank Quinn and the reviewer for
some essential comments and suggestions. 

\section{Main Results}

\begin{ssec}\label{centq}
Let $(A, m, \Delta , S, \epsilon, e)$ be a finite dimensional 
unimodular ribbon  Hopf algebra
over a field $k$ with an integral  $\Lambda \in A$
and a right integral $\lambda \in A^*$  such that $\lambda (\Lambda)=1$.
We define a linear map $\star \co A\otimes A\rightarrow A$ given by
$$
a\star b=\sum _b\lambda (S(a)b_{(1)})b_{(2)},
\;\;
\mbox{where}\;\; \Delta (b)=\sum_b b_{(1)}\otimes b_{(2)}.
$$
Let $Z(A)$ be the center of $A$ and let
$K(A)$ be the null space of the pairing on $Z(A)$ induced by
$\lambda$, i.e.
$$
K(A)=\{a\in Z(A)\mid\mbox{for any }\; b\in Z(A),\; \lambda(ab)=0\}.
$$
Then  $K(A)$ is an algebra ideal in $Z(A)$, and let
$\hat{Z}(A)=Z(A)/K(A)$ be the quotient algebra.
Given any $a\in Z(A)$, we will
denote by $[a]$ its equivalence class in $\hat{Z}(A)$.
Let also
$\hat{Z}^S(A)=\{ [a]\in \hat{Z}(A)\mid [S(a)]=[a]\}$ (this will be shown
to be well defined in \ref{well}).
\end{ssec}

\begin{lemma} \label{prop1} Let $A$ be a finite-dimensional unimodular
ribbon  Hopf algebra
over a field $k$ as above. Then
\begin{itemize}
\item[\rm(a)] $\star\,\co\, Z(A)\otimes Z(A)\rightarrow
Z(A)$ defines an associative product on $Z(A)$ with an identity $\Lambda$
and for any $a,b\in Z(A)$,
$
S(a\star b)=S(b)\star S(a);
$
\item[\rm(b)] $\star $ defines an associative and commutative product on
 $\hat{Z}(A)$.
\end{itemize}
\end{lemma}

\begin{ssec}\label{caa}
Let $C^n\subseteq A^{\otimes n}, \,n>1,$ be the centralizer of the action
of $A$ on $A^{\otimes n}$ given by
  the comultiplication, i.e.\ $a\in C^n$ iff for any $b\in
  A$, $ \Delta ^{n-1}(b)a=a\Delta ^{n-1}(b)$. Define also $C^1=Z(A)$.

$C^2$
contains the commutative subalgebra $C^2_Z$ generated by the elements of
the form $(a\otimes b)\Delta c$, where $a,b,c\in Z(A)$. Let
$
\mu \co C^2_Z\otimes C^2\rightarrow k,
$
be given by
$$
\mu( \sum_{i}a_i\otimes b_i,\sum_{j}c_j\otimes
d_j)=\sum_{i,j}\lambda(a_ic_j)\lambda(b_id_j),
$$
and let $\bar{\mu }\co C^2_Z\otimes C^2_Z\rightarrow k
$ be the corresponding restriction of $\mu$. Define
\beq
&&K^2_Z=\{x\in C^2_Z\mid \mu(x,y)=0\;\mbox{ for
any}\; y\in  C^2\}\; \mbox{ and }\\
&&\overline{K}^2_Z= \{x\in C^2_Z\mid \mu(x,y)=0\;\mbox{ for
any}\; y\in  C^2_Z\}.
\eeq
Obviously $K^2_Z$ and $\overline{K}^2_Z$ are ideals in
$C^2_Z$ and $K^2_Z\subset \overline{K}^2_Z$.
This induces a surjective homomorphism
$$
\pi_Z\co \,C^2_Z/K^2_Z\rightarrow C^2_Z/\overline{K}^2_Z.
$$
\end{ssec}
Define
$
\delta \co Z(A)\otimes Z(A)\rightarrow C_Z^2$ as $\delta(w,z)=z\otimes
w-(1\otimes w)\Delta (z)$.
\begin{prop} \label{wd} $\delta$ factors through a well defined map
$\hat{\delta}\co\hat{Z}(A)\otimes \hat{Z}(A)\rightarrow C^2_Z/K^2_Z$.
\end{prop}

\begin{ssec}
Let $ \T\subset \T_Z\subset\hat{Z}^S(A)$ be
$
\T=\{[z]\in \hat{Z}^S(A)\mid \hat{\delta}([z],[z])=0\}$ and $
\T_Z=\{[z]\in \hat{Z}^S(A)\mid \pi_Z\cdot\hat{\delta}([z],[z])=0\}$.
 Observe that
$[z]\in \T_Z$
if and only if for any
$a,b,c\in Z(A)$,
$
\lambda(zc(bz\star a))=
\lambda(zc(b\star (za)))
$.
Hence
\end{ssec}
\begin{prop} \label{prop2}
    $[z]\in \T_Z$ if and only if
for any
$[a],[b]\in \hat{Z}(A)$,
$
[z(a\star zb)]=[z(az\star b)].
$
\end{prop}

\begin{ssec}\label{Jdef}
    Let
 $J\co Z(A)\rightarrow Z(A)$, be defined as
$$
J(z)=(\lambda \otimes  1) (z\otimes
1)R^{21}R=\sum_{i,j}\lambda(z\beta_i\alpha_j)\alpha_i\beta_j.
$$
This operator is related to the image of one of the generators, ${\cal
S}$, in the action of the torus group
on $Z(A)$ (see \cite{K:mcg}, (2.55))
and it is
essential in understanding when the invariant of the 4--thickening
reduces to an invariant of the boundary.  Let
 $Z_{\star}(A)$
denote the algebra which has $Z(A)$ as a vector space  and the $\star $
product structure. Then
\end{ssec}
\begin{prop} \label{J}\begin{itemize}
 \item[{\rm(a)}] $J\co Z_{\star}(A)\rightarrow Z(A)$ is an
algebra homomorphism, i.e.\ for any $a,b\in Z(A)$,
$
J(a\star b)=J(a)J(b).
$
\item[\rm(b)] $J^2(a)=S(a)\star J(1)$;
\item[\rm(c)] $J$ factors through an algebra homomorphism map
$\hat{J}\co\hat{Z}_{\star}(A)\rightarrow
\hat{Z}(A)$, and  maps $\hat{Z}^S(A)$ into itself.
\end{itemize}
\end{prop}

Observe that, if $J(1)=\gamma \Lambda$, where $\gamma\in k$ is a unit,
\ref{J} (b) and the fact that on the center of a ribbon algebra 
$S^2$ acts as the identity, imply
that $J$ is bijective with an inverse $J^{-1}=\gamma ^{-1}(S\circ J)$.
Then from
\ref{J} (a) and \ref{prop1} (a) one obtains
\beq
&J(ab)&=J(J\circ J^{-1}(a)J\circ J^{-1}(b))=J^2(J^{-1}(b)\star
J^{-1}(a))\\
&&=\gamma ^{-1}S(S\circ J(b)\star S\circ J(a))=\gamma ^{-1}J(a)\star J(b).
\eeq
Therefore we have proved the following:
\begin{corr}\label{corrJ} If $J(1)=\gamma\Lambda$, where $\gamma\in k$ is a
unit, then
  $\gamma ^{-1}J\co Z(A)\rightarrow Z_{\star}(A)$ is an algebra
isomorphism. In particular, the algebra $Z_{\star}(A)$ is commutative.
\end{corr}

\begin{dfn}
    A quasitriangular unimodular ribbon Hopf algebra for which 
    $J(1) =\gamma \Lambda$,
 where $\gamma\in k$ is a unit, will be called  $\Lambda$--{\em
 factorizable}.\footnote{A quasitriangular Hopf algebra is called
 factorizable if $\bar{J}\co A^*\rightarrow A$, given by
 $\bar{J}(f)=(f\otimes 1)(R^{21}R)$ is bijective.}
\end{dfn}

\begin{lemma} \label{Jinv}  $\T$ is a commutative monoid with respect to
the usual
and the $\star$--product on $\hat{Z}^S(A)$. Moreover
$\hat{J}$ sends $\T$ into itself.
\end{lemma}

We observe that proposition \ref{J} implies that
when the algebra is $\Lambda$--factorizable, $\hat{J}\co\T\rightarrow \T$
is a bijection whose square is a multiple of the identity.

\begin{ssec} 
    Let $M$ be an orientable 4--dimensional manifold
which possesses a
  decomposition as a handlebody with 0--, 1-- and 2--handles.
  We remind that an $n$--handle is a product $D^n \times D^{4-n}$ and
the choice of  radial 
coordinates in $D^{4-n}$ gives a description of the product as the mapping 
cylinder of a projection $D^n \times S^{{3-n}} \rightarrow D^n$. 
Then $D^n \times \{0\}$ is called the core, 
$S^{n-1}\times \{0\}$ is called the attaching sphere and $\{0\}\times 
S^{{3-n}}$ is called the belt sphere of the handle.
When another handle is 
attached on top of this one the intersection of the attaching map with 
the handle lies in $D^n \times S^{3-n}$ and 
using the mapping cylinder coordinates
the core of the upper handle 
 can be extended in the lower handle. This 
extends the upper cores to a disk whose boundary lies on the lower 
cores. The union 
of these extended cores forms a 
2--dimensional CW complex which  will be called the {\em spine} of the 
handlebody.  The mapping cylinder contractions also combine to give a 
standard deformation retraction of the handlebody to the spine.

A pair of $(n+1)$--handle and an $n$--handle is called a {\em cancelling pair} if 
  the attaching sphere of the $(n+1)$--handle intersects the belt sphere 
     of the $n$--handle in a single point.
     
Then a {\em  4--thickening} $M$ of a 2--dimensional CW complex $P$, denoted with
$(M,P)$, is an orientable 4--dimensional manifold
together with a
  decomposition as a handlebody with 0--, 1-- and 2--handles
and an identification (as CW complexes)
of the spine of the handlebody structure with
$P$ through an embedding $\iota_{M,P}\co P\rightarrow M$. 
In particular, $\iota_{M,P}$ induces  isomorphism on homology.
We will restrict ourselves to  4--thickenings with a single 0--handle.
A 2--deformation of such 4--thickenings is given by a sequence of the 
following handle moves:
\begin{itemize}
    \item[(a)]  creation or 
    cancellation of a cancelling 1--2 handle pair;   
    \item[(b)] changing the attaching maps of the 1-- and 2-- handles by 
    isotopy.
    \end{itemize}
 Observe that these moves induce a 2--deformation on the spine.
 
 The word 4--thickening is supposed to stress not only the fact that
a spine has been fixed, but also that we have weakened the equivalence
relations on the objects with respect to 4--manifolds.\footnote{While 
changing the attaching map of a 2--handle by isotopy  is 
equivalent to the creation and  cancellation of cancelling 2--3 
handle pairs, isotoping   the attaching map of a 3--handle  is 
not a 2--deformation.}
\end{ssec}

\begin{ssec} The monoid $\T$ will be shown to
correspond to invariants under 2--handle
slides. An invariance under 2--deformations requires in addition
invariance under 1--2 handle cancellations, and the center elements
which lead to such invariants form the following subset of $\T$:
$$
\T^4=\{[z]\in \T\mid \mbox{there exits } [w]\in \hat{Z}^S(A) \mbox{ and }
[zw]=[\Lambda]\}.
$$
Let also 
$\T^3=\{[z]\in \T^4 \mid  [zJ(z)]=X_z[\Lambda] \mbox{ for some unit }
X_z\in k\}$ and
\begin{gather*}
\T^2=\{[z]\in \T^4\mid [z]=[z_1J(z_2)] \mbox{ and }\hat{\delta}([z_1],[z_2])=0
\hspace{1.6cm}\\\hspace{6cm}\mbox{ for some } [z_1],[z_2]\in 
\hat{Z}^S(A)  \}.\end{gather*}
\end{ssec}

  \begin{theo} \label{theo1} Given any
$[z]\in \T^4$ and $[w]\in \hat{Z}^S(A)$ such that $[zw]=[\Lambda]$,
there exists a HKR--type invariant of 4--thickenings under
2--deformations, denoted with $\Z_{[z]} (M)$, such that
 $$
   \Z_{[ z]} (S^{2}\times D^{2})=\lambda (z)\; \mbox{and}\;
   \Z_{[ z]} (S^{1}\times D^{3})=\epsilon(w ).
   $$
\end{theo}

Obviously for any finite dimensional unimodular ribbon Hopf algebra 
$A$, the elements $[1],[\Lambda]\in \T^4$. The choice $[z]=[\Lambda]$ brings to
the trivial invariant which is 1 for any $M$.
On another hand $[z]=[1]$ gives the Hennings invariant (in the 3--manifold
case):

\begin{corr} \label{corr1} Any  finite-dimensional unimodular  ribbon
       Hopf algebra   $A$ over a field $k$,
    determines an invariant $\Z_ A$ of 4--thickenings under
   2--deform\-at\-ions, such that
   $$
   \Z_ A(S^{2}\times D^{2})=\lambda (1),\; \mbox{and}\;
   \Z_ A(S^{1}\times D^{3})=\epsilon(\Lambda ),
   $$
   In particular,  $\Z_ A(S^{2}\times D^{2})\neq 0$ if and only if $A$ is
   cosemisimple ($A^{*}$ is semisimple), and $\Z_ A(S^{1}\times D^{3})\neq
0$ if and
   only if $A$ is semisimple.
   \end{corr}

   Given a 4--manifold $M$, let $w_2(M)\in H^2(M;Z/2)$ denote the
   second Whitney class of $M$.

   \begin{lemma} \label{lemma2}
Let
$P$ be a 2--dimensional CW complex and $(M_1, P)$, $(M_2, P)$ be two
4--thickenings of $P$ such that
$\iota_{M_1,P}^*(w_2(M_1))=\iota_{M_2,P}^*(w_2(M_2))$.
If  $[z]\in \T^2$ then
$\Z_{[ z]}(M_1)=\Z_{[ z]}(M_2) .$
   \end{lemma}

\begin{corr}\label{corr2}
Let  $A$  be a triangular Hopf algebra
and let $[z]\in \T^4$. If
 $(M_1,P_1)$ and $(M_2,P_2)$ are two 4--thickenings such that
 $P_1$ and $P_2$ are related by a 2--deformation, then 
$\Z_{[z]}(M_1)=\Z_{[z]}(M_2)$.
   \end{corr}
Hence, if $A$ is a triangular Hopf
algebra any $[z]\in\T^4$ defines an invariant
of 2--complexes under 2--deformations, and this invariant  is
 denoted by $\Z_{[z]}^2(P)$. 
Then it is natural to expect that for triangular algebras
$\T^4=\T^2$. Actually, in this case for any $z\in Z(A)$,
 $J(z)=\lambda(z)1$. In particular,
$$\T^2=\{[z]\in \T^4\mid \mbox{ there exists } [w]\in 
\hat{Z}^S(A)\mbox{ with }\hat{\delta}([z],[w])=0 \mbox{ and } 
\lambda(w)\neq 0\}.
$$
And since for any $z\in Z(A)$, $\delta(z,1)=0$,  it follows that
if $A$ is  triangular and cosemisimple (i.e.\ $\lambda(1)\neq 0$) then 
$\T^2=\T^4$. We don't know if this is true for any triangular algebra.

\begin{ssec} Let $M$ be a 4--thickening represented with a Kirby diagram $L$ (see
section 5) and  let $\sigma_{+}$,
$\sigma_{-}$ and $\sigma _0$
be the
numbers of positive,
negative and zero  eigenvalues of the linking matrix of $L$.
\end{ssec}

\begin{corr} \label{3inv}
If  $[z]\in \T^3$  then
$C_{+ }=\Z_{[ z]} (CP^2)$
and $C_{- }=\Z_{[ z]} (\overline{CP^2})$ are units in $k$. Moreover,
  if $M$ is a 4--thickening with $n$ 1--handles,  then
$$
C_+^{n-\sigma_{+}}C_-^{n-\sigma_{-}}\Z_{[ z]} ( M)
$$
only depends on the boundary $\partial M$ of $M$ and is denoted by
$
\Z_{[ z]}^{\partial}(\partial M)
$.
\end{corr}

   \section{Basic facts about Hopf algebras}

Here, we introduce some notations assuming that the reader is familiar with the
axioms of a Hopf algebra. A possible reference about Hopf algebras
is \cite{Sw}. Let $(A, m, \Delta , S, \epsilon, e)$ be a Hopf algebra
over a field $k$, where:
\begin{eqnarray*}
      && m\co A\otimes A \rightarrow A \quad \mbox{ multiplication map}\\
       &&\Delta\co A\rightarrow  A\otimes A \quad \mbox{ comultiplication map}\\
       &&S\co A \rightarrow A ^{opp}\quad \mbox{ antipode}\\
      && \epsilon\co A\rightarrow k \quad \mbox{ counit}\\
      && e\co k\rightarrow A  \quad \mbox{ unit}
   \end{eqnarray*}
Note also that there are natural isomorphisms; $k\otimes A\rightarrow
A$ and $A\otimes k \rightarrow A$ which we will often omit,
identifying $A\otimes k$ and $k\otimes A$ with $A$.

\begin{ssec} \label{propH}
The maps above need to satisfy a list of compatibility
conditions, out of which we only mention the following:
\begin{itemize}
      \item[(a)]$\Delta (\Delta \otimes 1)=\Delta (1\otimes \Delta
      )\co A\rightarrow A\otimes A\otimes
      A\quad\mbox{ (coassociativity)}$,
     \item[(b)]$\Delta  m =(m\otimes m)(1\otimes T\otimes 1)(\Delta \otimes
      \Delta) \co A\otimes A\rightarrow A\otimes A,\quad  \Delta(1)=1\otimes
      1$,
      \item[(c)]$ m (S\otimes 1)\Delta =m (1\otimes S)\Delta =e \epsilon
      \co A\rightarrow A$,
    \end{itemize}
    where 1 denotes both the identity element $e(1_{k})$ in $A$ and the identity
    map $A\rightarrow A$, and $T\co A\otimes A\rightarrow A\otimes A$ is the
     transposition map $a\otimes b\rightarrow b\otimes a$.
An easy consequence of the definition of the antipode is that
\begin{itemize}
\item[(d)] $T\circ (S\otimes S)\Delta(a)=\Delta (S(a))$.
\end{itemize}
     Let
     $
     \Delta ^n =(\Delta\otimes 1^{\otimes (n-1)})(\Delta \otimes 1^{\otimes
     (n-2)}) \ldots \Delta \co  A\rightarrow A^{\otimes (n+1)}.
     $
We use Sweedler's notation
     $
     \Delta ^{(n-1)} (a)=\sum_{a}a_{(1)}\otimes a_{(2)}\otimes \ldots
     a_{(n-1)}\otimes a_{(n)}
     $. Then (d) implies that
\begin{itemize}
\item[(e)]    $\Delta ^{n-1}(S(a))= \sum_{a}S(a_{(n)})\otimes S(a_{(n-1)})
     \ldots \otimes S(a_{(1)}).$
     \end{itemize}
     \end{ssec}

   \begin{ssec}   An element $\lambda _L\in A^{*}$ is called a {\em left
integral for} $A^*$ if
      $$
      (f\otimes \lambda _L)\Delta(a)=\lambda _L(a) f(1), \; \mbox{for any
$a\in A$ and
      $f\in A^*$.}
     $$
      An element $\lambda _R\in A^{*}$ is called a {\em right integral
      for} $A^*$ if
      $$
      (\lambda _R\otimes f)\Delta(a)=\lambda _R(a) f(1), \; \mbox{for any
$a\in A$ and
      $f\in A^*$.}
      $$

     When $A$ is
     finite-dimensional, the Hopf algebra isomorphism $A\simeq A^{**}$
implies that
     one can define a left (right) integral
     for $A$ as an element $\Lambda \in A$, such that $a.\Lambda
     =\epsilon (a)\Lambda $ ($\Lambda a
     =\epsilon (a)\Lambda $) for any $a\in A$.
    \end{ssec}

   \begin{ssec} \label{print}
       The following results (\cite{Sw, R:trace, R:ant}) concern
the existence of
      integrals when $A$ is a finite-dimensional Hopf algebra over
      a field $k$.
      \begin{itemize}
	\item[(a)] The subspaces $\int _L^*, \int _R^*\subset A^*$
	of left (right) integrals for $A^*$ and the subspaces
            $\int _L, \int _R \subset A$
	    of left (right) integrals for $A$
	  are one dimensional;
	\item[(b)] The antipode map is bijective;
	\item[(c)] For any nonzero $\lambda \in \int _R^*$ there exists
	       $\Lambda \in \int_L $ such that
           	$$\lambda (\Lambda )=\lambda(S(\Lambda))=1;$$
	\item[(d)] Given any nonzero $\lambda\in \int_R^*$
          	the map $\Phi \co A\rightarrow A^*$ given by
	 $\Phi(a)(b)=\lambda (a b)$ is a bijection;
	\end{itemize}
    \end{ssec}

   \begin{ssec}\label{unimodular}
	Note that, if A is a finite-dimensional Hopf algebra and $\Lambda\in
	\int _{R}$, then $S(\Lambda),\, S^{-1}(\Lambda)\in \int_{L}$.
	Moreover, if $\lambda \in \int _R^*$, then $\lambda \circ S,
	\lambda\circ S^{-1}\in \int _L^*$.
	$A$   is called {\em unimodular} if $\int_R=\int_L$
	and if $A$ is unimodular then for any $\lambda\in \int _R^*$,
$\lambda(a b)=\lambda (S^2(b)a).$
	\end{ssec}

   \begin{ssec}\label{quasitriangular}
   A {\em quasitriangular} Hopf algebra is a Hopf algebra
   $A$ endowed with
   invertible element $R=\sum_{i}\alpha_{i}\otimes \beta_{i}\in A\otimes A$
   such that
   \begin{itemize}
   \item[(a)] $T\circ \Delta (a) = R\Delta (a) R^{-1}\;\mbox{for any $a\in
   A$}$;
   \item[(b)] $(\Delta \otimes 1)R=R^{13}R^{23}$;
    \item[(c)] $(1 \otimes \Delta)R=R^{13}R^{12}$,
\end{itemize}
      where as usual $R^{(kl)}\in A^{\otimes n}$ indicates the image
      of $R$ under the injective homomorphism of the group of
      invertible elements in $
      A\otimes A$ into the group of invertible elements of $A^{\otimes n}$
      where the first factor is mapped into $k$-th position and the
      second into $l$-th position.

If $(A,R)$ is a quasitriangular Hopf algebra, the following relations
hold:
\begin{itemize}
      \item[(d)] $R^{(12)}R^{(13)}R^{(23)}=R^{(23)}R^{(13)}R^{(12)}$;
      \item[(e)] $(S\otimes 1)R=(1\otimes S^{-1})R=R^{-1}$, and
      $(S\otimes S)R=R$;
      \item[(f)] $(\epsilon\otimes 1)R=(1\otimes \epsilon )R=1$;
      \item[(g)] Let $u=\sum_{i}S(\beta_{i})\alpha_i$, then $u$ is
      invertible and $S^2(a)=uau^{-1}$, moreover,
      $$
      \Delta (u)=(u\otimes u)(R^{(21)}R)^{-1}.
      $$
      \end{itemize}
 \end{ssec}

\begin{ssec} A quasitriangular Hopf algebra is called {\em
triangular} if $R^{-1}=R^{(21)}=\sum_{i}\beta_{i}\otimes \alpha_{i}$.
In this case
  $u$
  is a group-like element, i.e.\
  $\Delta  (u)=u\otimes u$, which, in the terminology below,
  implies that any triangular Hopf algebra is ribbon with ribbon
  element $u$.

 A  Hopf algebra $A$ is called {\em cocommutative} if it possesses
 triangular structure with      $
     R=1\otimes 1
     $, i.e.\ if $T\circ \Delta =\Delta$.
\end{ssec}

\begin{ssec}\label{ribbon} A quasitriangular Hopf algebra $A$ is called
{\em ribbon}
if it is endowed with a  grouplike element $g\in A$ such that
$S^2(a)=gag^{-1}$, called the special grouplike element of $A$
(grouplike means that $g$ is invertible and $\Delta g=g\otimes g$).
It can be shown (see for example \cite{RT, KR}) that if $A$ is ribbon,
$$ \theta=gu^{-1}=u^{-1}g=
\sum _i \alpha_{i}g^{-1} \beta_{i}=\sum _i\beta_{i} g\alpha_{i}
$$
is a central  element in $A$ such that
\begin{itemize}
\item[(a)] $S(\theta)=\theta$;
\item[(b)] $\theta$ is invertible
with  inverse
$
\theta^{-1}=\sum _i\alpha _iS(\beta_i) g=\sum_iS(\beta_i)\alpha_i g^{-1};
$
\item[(c)] $\Delta (\theta)=(\theta\otimes \theta)(R^{(21)}R)^{-1}.
      $
      \end{itemize}
$\theta$ is called the
{\em ribbon element of $A$}.

A {\em trace function } on $A$ is an element $f\in A^{*}$ such that,
for any $a,b\in A$,
$f( ab)=f(ba)$ and $f(a)=f(S(a))$.
In a finite dimensional unimodular
ribbon Hopf algebra there is a bijection between the set of
$S$--invariant central elements in $A$ and the space of trace functions
on $A$ given by  $z\rightarrow \lambda _{zg}$, where
 $\lambda _{zg}(a)=\lambda (zga )$ (\cite{Hennings,R:trace}).
\end{ssec}

\section{The center of a unimodular finite dimensional ribbon Hopf
algebra}

 In the rest of the paper, unless specified otherwise,
  $(A, m, \Delta , S, \epsilon, e)$ will be a unimodular  Hopf algebra
over a field $k$ with an integral  $\Lambda \in A$,
 a right integral $\lambda \in A^*$ and a left integral $\lambda
 ^S=\lambda\circ S$, such that $\lambda
(\Lambda)=\lambda ^S(\Lambda)=1$. Moreover, we assume that $A$
carries a ribbon structure
 given by an $R$--matrix
$R=\sum_i\alpha_i\otimes \beta_i$ and a group like element $g$ such
that $gag^{-1}=S^2(a)$ for any $a\in A$.
 Many of the statements here can be easily
illustrated using the diagrammatic language in the later chapters, but
because of their purely algebraic significance we decided that it is
better to prove them in a self-contained way.

\begin{ssec} \label{gencn}
{\bf Generating elements in $C^n$}
\begin{itemize}
\item[(i)] The first way to
generate elements in  $C^n$, is by ``going up'', i.e.\ by
applying some of the following embeddings on $C^{n-1}$:
\beq
&&\eta ^{(n-1)}_r \co C^{n-1}\rightarrow C^n,\; a\rightarrow 1\otimes
a;\\
&&\eta ^{(n-1)} _l\co C^{n-1}\rightarrow C^n,\; a\rightarrow a\otimes
1;\\
&&1^{\otimes (i-1)}\otimes \Delta \otimes 1^{\otimes
(n-i-1)} \co C^{n-1}\rightarrow C^n,\; i=1, \ldots, n-1.
\eeq
The subalgebra of $C^n$ generated inductively in this way,
 starting with $C^1=Z(A)$, will be denoted with $C^n_Z$.

\item[(ii)]
The second way to generate new elements in $C^n$ is through the action of
the braid group on $C^n$ as follows. If $B_n$ is the
braid group on $n$ strings and $q_n\co B_n\rightarrow {\bf S_n}$ is its
homomorphism onto the symmetric group ${\bf S_n}$, let $I_n=q_n^{-1}(id)$.
The relation  \ref{quasitriangular} (d) implies that one
 can define a representation of $\phi \co B_n\rightarrow End(A^{\otimes n})$ by
defining the image of the generator which interchanges the $i$-th and
the $(i+1)$-st strings to be
$$
\phi (\sigma_{i,i+1})=1^{\otimes (i-1)}\otimes (T\circ R)\otimes 1^{\otimes
(n-i-1)},
$$
where we first multiply the corresponding element in $A^{\otimes n}$
on the left with $1^{\otimes (i-1)}\otimes R\otimes 1^{\otimes
(n-i-1)}$ and then apply the permutation.
Suppose that $s,s'\in B_n$ are such that $q_n(s)=q_n(s')^{-1}$.
Then the condition \ref{quasitriangular} (a) implies that given any
$a\in C^n$,
$\phi (s)\circ a\circ \phi (s')$ act on $A^{\otimes n}$ by multiplication with
an element in  $C^n$. We write this fact as
$
\phi (s)\circ C^n\circ \phi (s')\subset C^n
$. For example, if $\sum_i c_i\otimes d_i\in C^2$ then
$
\sum_{i,k,j} \beta_k d_i\alpha_j\otimes \alpha_k c_i\beta_j\in C^2$.
The statement implies in particular that $\phi (I_n)\subset C^n$.
\footnote{Using \ref{ribbon} (c) one can show that actually
$\phi (I_n)\subset C^n_Z$.}

\item[(iii)] The third way to obtain elements in $C^n$ is by ``going
down'', i.e.\ by applying the
integrals to the elements in $C^{n+k}$:
\end{itemize}
\end{ssec}
\begin{prop}\label{cm}
    Let ${\cal L}_{n+1}\co A^{\otimes (n+1)}\rightarrow A^{\otimes n}$
    be the map which applies $\lambda $ on the leftmost factor in 
    $A^{\otimes (n+1)}$ and
    let ${\cal R}_{n+1}\co A^{\otimes (n+1)}\rightarrow A^{\otimes n}$
    be the map which applies $\lambda ^S$ on the rightmost factor in 
    $A^{\otimes (n+1)}$.
    Then ${\cal L}_{n+1}$ and ${\cal R}_{n+1}$ map $C^{n+1}$ into
    $C^n$.
\end{prop}
\begin{proof}
The proof is standard, but for completeness we will show the first
part of the statement and the second is analogous.
Given any $\sum_i a_i\otimes b_i \in C^{ n+1}$,
where $a_i\in A,\; b_{i}\in A^{\otimes n}$ and any $c\in A$,
\begin{eqnarray*}
&\sum_i\lambda (a_i)b_i\Delta^{n-1}(c)&=\sum_{i,c}\lambda
(a_ic_{(2)}S^{-1}(c_{(1)}))b_i\Delta^{n-1}(c_{(3)})\\
&&=\sum_{i,c}\lambda
(c_{(2)}a_iS^{-1}(c_{(1)}))\Delta^{n-1}(c_{(3)})b_i\\
&&=\sum_{i,c}\lambda
(S(c_{(1)})c_{(2)}a_i)\Delta^{n-1}(c_{(3)})b_i=\sum_i\lambda
(a_i)\Delta^{n-1}(c)b_i,
\end{eqnarray*}
hence $\sum_i\lambda (a_i)b_i\in C^n$. \end{proof}

By induction the last proposition implies that
for any $0\leq k< l\leq n$
$$\lambda^{\otimes k}\otimes 1^{\otimes (l-k)}\otimes
(\lambda^S)^{\otimes (n-l)}\co C^n\rightarrow
C^{l-k}.
$$

\begin{prop} \label{Sinv} For any  $a\in C^n$ and any partition
$n'+n''=n$,
 $
 \lambda ^{\otimes n}(a)=(\lambda ^{\otimes n'}\otimes
 (\lambda ^S) ^{\otimes n''}) (a)
 $. In particular, $\lambda (a)=\lambda ^S(a)$ for any $a\in Z(A)$.
 \end{prop}
 \begin{proof}
 First we will prove the statement for $n=1$.
 Suppose that $a\in Z(A)$. Then, using \ref{ribbon}, it follows that
\begin{eqnarray*}
&\lambda (a)&=
\sum_{i,j}\lambda (a\beta_j\beta_i S^{-1}(\alpha _i)\alpha_j)=
\sum_{i,j}\lambda (a\beta_i S^{-1}(\alpha
_i)\alpha_jS^{-2}(\beta_j))\\
&&=\sum_{i,j}\lambda (gagg^{-1}\beta_i S^{-1}(\alpha
_i)\alpha_jg^{-1}\beta_j)=
\lambda (gagS^{-1}(\theta^{-1})\theta)\\
&&=\lambda (gag)=\lambda(S(a)).
\end{eqnarray*}
Let now  $a\in C^{n}$, $n>1$. If $n''=0$, the statement is
trivial. Suppose then that it is true for some $n''\geq 0$.
 Then proposition \ref{cm} implies that $(\lambda ^{\otimes (n'-1)}\otimes
 1\otimes ( \lambda ^S)^{\otimes n''})(a)\in
 Z(A)$ and hence the statement with $n''+1$ follows from the one for
 $n''$ and from the statement with $n=1$.\end{proof}
This proposition implies that if $a\in K(A)$ then for any $b\in Z(A)$,
$
\lambda (bS(a))=\lambda (S^2(a) S(b))=\lambda (S(b)a)=0
$, i.e.\ $S(a)\in K(A)$. Hence

\begin{corr}\label{well} The algebra $\hat{Z}^S(A)$ in
\ref{centq} is well defined.
\end{corr}

\begin{ssec} \label{proofp1} {\bf Proof of lemma \ref{prop1}}\qua
First observe  that proposition \ref{cm} implies that,
for any $a,b\in Z(A)$, $a\star b\in Z(A)$.
To see the associativity of the product, let
 $a,b,c\in Z(A)$. Then
\beq
&(a\star b)\star c&=\sum_{c,b}\lambda (S(a) b_{(1)})\lambda
(S(b_{(2)})c_{(1)})c_{(2)}\\
&&=\sum_{c,b}\lambda (S(a) b_{(1)}S(b_{(2)})c_{(2)})\lambda
(S(b_{(3)})c_{(1)})c_{(3)}\\
&&=\sum_{c}\lambda(S(a)c_{(2)})\lambda(S(b)c_{(1)})c_{(3)}=a\star(b\star c).
\eeq
To complete the proof of
\ref{prop1}(a) we observe that for any $a,b\in Z(A)$,
 \beq
&S(a\star  b)&=\sum_{b}\lambda (S(a)b_{(1)})S(b_{(2)})=
\sum_{S(a),b}\lambda(S(a)_{(1)}b_{(1)})S(a)_{(2)}b_{(2)}S(b_{(3)})\\
&&=\sum_{S(a)}\lambda (S(a)_{(1)}b)S(a)_{(2)}=S(b)\star  S(a),
\eeq
which together with the definition of $\Lambda $ implies that
$
a=\Lambda \star a=a\star \Lambda
$.
This completes the proof of proposition \ref{prop1} (a). Now, for any $a,b,c
\in Z(A)$, define
$$
\sigma (a,b,c)=\lambda (S(a)(b\star c)).
$$
Then \ref{prop1} (b)  follows from \ref{Sinv} and the following proposition.
\end{ssec}
\begin{prop}\label{strco}
If $(a',b',c')$ is any permutation of $(a,b,c)$ or\hfill\break
$(S(a),S(b),S(c))$, then $$\sigma (a',b',c')=\sigma (a,b,c).$$
\end{prop}
\begin{proof}
First we observe that \ref{prop1} (a) and \ref{Sinv} imply that
$$
\sigma(a,b,c)=\sigma(S(a),S(c),S(b))
.$$
Hence, it is enough to show that $\sigma (a',b',c')=\sigma (a,b,c)$
where $(a',b',c')$ is one of the two permutations $(b,a,c)$ or
$(c,b,a)$.  Now we claim  that $\lambda (a(S(b)\star
c))=\lambda (b(S(a)\star c))$ which would
imply that $\sigma (a,b,c)=\sigma (b,a,c)$.
To see this, let $\sum_i\gamma_i\otimes \delta _i=R^{-1}$. Then
\begin{eqnarray*}
&&\lambda (b(S(a)\star c))=\sum_{c}\lambda (ac_{(1)})\lambda (bc_{(2)})=
\sum_{c,i,j}\lambda (a\gamma_i\alpha_j c_{(1)})\lambda (b\delta_i\beta_jc_{(2)})
\\
&&=\sum_{c,i,j}\lambda (a\gamma_i c_{(2)}\alpha_j)\lambda (b\delta
_ic_{(1)}\beta_j)=
\sum_{c,i,j}\lambda (a c_{(2)}\alpha_jS^{-2}(\gamma_i))\lambda
(bc_{(1)}\beta_jS^{-2}(\delta
_i))\\
&&=\sum_{c}\lambda (bc_{(1)})\lambda (ac_{(2)})=\lambda(a(S(b)\star c)).
\end{eqnarray*}
We  complete the proof of the proposition as follows:
$$
\sigma(a,b,c)=\sigma(S(a),S(c),S(b))=\sigma(S(c),S(a),S(b))=\sigma(c,b,a).
$$\end{proof}

\begin{ssec} {\bf Proof of proposition \ref{wd}}\qua
It is enough to show that for any $z\in K(A)$ and  any
$\sum_i a_i\otimes b_i \in C^2$, the following three statements hold:
\begin{itemize}
\item[(a)]$\sum_i \lambda( a_i)\lambda(z b_i)=0$,
\item[(b)]$\sum_{z,i} \lambda( z_{(1)}a_i)\lambda(z_{(2)} b_i)=0$,
\item[(c)]$\sum_i \lambda( za_i)\lambda( b_i)=0$.
\end{itemize}
(a) and (c) follow directly from  \ref{Sinv} and \ref{cm}.
On another hand to show (b), using  \ref{Sinv} and the fact that
$z=z\star \Lambda$, we obtain
\begin{eqnarray*}
&\sum_{z,i} \lambda( z_{(1)}a_i)\lambda(z_{(2)} b_i)&
=\sum_{\Lambda,i} \lambda(S(z) \Lambda_{(1)})
\lambda( \Lambda_{(2)}a_i)\lambda(\Lambda_{(3)}
b_i)\\
&&=\sum_{\Lambda} \lambda(z(\sum_i
\lambda^S(
\Lambda_{(2)}a_i)\lambda^S(\Lambda_{(3)}b_i)S(\Lambda_{(1)})))=0.
\end{eqnarray*}
\end{ssec}

\begin{ssec} \label{exJ}{\bf Proof of proposition \ref{J}}\qua
Observe that $J$ actually maps the center into itself since
from \ref{ribbon} it follows that
$$J(z)=\sum_{\theta}\lambda (z\theta \theta^{-1}_{(1)})\theta \theta^{-1}_{(2)}
=\theta((S(z)\theta)\star
(\theta^{-1})).
$$
This expression also implies (together with \ref{prop1} (b) )
that $J$ factors through a map $\hat{Z}(A)\rightarrow
\hat{Z}(A)$. Now we can complete the proof of \ref{J} (c).
Let $[a]\in \hat{Z}^S(A)$. Then
using the fact that $S(\theta)=\theta$ and \ref{prop1} (a) and (b)
we obtain that
$$
[S(J(a))-J(a)]=[\theta(((S(a)-a)\theta)\star (\theta^{-1}))]=0.
$$
Hence $[J(a)]\in \hat{Z}^S(A)$.

It is left to show \ref{J} (a) and (b).
\begin{itemize}
\item[(a)] Let $J'=J\circ S$. Then
(a) is equivalent to show that $J'\co Z_{\star}(A)\rightarrow Z(A)$
is an algebra isomorphism, i.e.\ for any $a,b\in Z(A)$,
$J'(a\star b)=J'(a)J'(b)$. From \ref {quasitriangular} (b) and (c) it
follows that
\beq
&J'(a)J'(b)&=\sum_{i,j}\lambda(S(a)\beta_i\alpha_j)
\alpha_i J'(b)\beta_j\\
&&=\sum_{i,j,i',j'}\lambda(S(a)\beta_i\alpha_j)
\lambda(S(b)\beta_{i'}\alpha_{j'})\alpha_i\alpha_{i'}\beta_{j'}\beta_j\\
&& =\sum_{i,j,\alpha
_j,\beta_i}\lambda(S(a)\beta_{i,(2)}\alpha_{j,(2)})
\lambda(S(b)\beta_{i,(1)}\alpha_{j,(1)})\alpha_i\beta_{j}\\
&&
=\sum_{i,j,b}\lambda(S(a)b_{(1)})\lambda(S(b_{(2)})\beta_i\alpha_j)\alpha_i
\beta_j=
J'(a\star b).
\eeq
\item[(b)] From \ref {quasitriangular} (b) and (c) it follows that
\beq
&S(a)\star J(1)&=\sum_{i,j,k,l}\lambda(\beta_i\beta_k\alpha_l\alpha_j)
\lambda(a\alpha_i\beta_j)\alpha_k\beta_l\\
&&=\lambda(\alpha_j\beta_i\beta_k\alpha_l)\lambda(a\beta_j\alpha_i)\alpha_k
\beta_l=
J^2(a).
\eeq
\end{itemize}
\end{ssec}

\begin{ssec} {\bf Proof of lemma \ref{Jinv}}\qua
It is obvious that $\T$ is a monoid under the usual
multiplication in $\hat{Z}^S(A)$.

First we will show that if
$\hat{\delta}([z],[z])=0$ and $\hat{\delta}([w],[w])=0$, then \\
$\hat{\delta}([z\star w],[z\star w])=0$.
This is equivalent to say that  for any $\sum_k a_k\otimes
b_k\in C^2$,
$$
\lambda (x(z\star w))=
\sum_{k}\lambda ((z\star w)a_k)\lambda ((z\star w)b_k),
$$
where $x=\sum_{k,w}\lambda (S(z) w_{(1)})\lambda(w_{(2)}a_k)w_{(3)}b_k\in
Z(A)$. From \ref{strco} it
follows that the left hand side is actually equal to
$\sigma(S(x),z,w)=\sigma(S(w),S(z),x)$. Hence
\beq
&l.h.s&=  \sum_{k,w}\lambda (S(z) w_{(1)})\lambda(w_{(2)}a_k)\lambda
(z w_{(3)}b_{k,(1)})\lambda(w w_{(4)}b_{k,(2)})\\
&&=\sum_{k,w,i,j}\lambda (S(z) w_{(1)})\lambda(\alpha_i
w_{(2)}a_kS(\alpha _j))\lambda
(z \beta_i w_{(3)}b_{k,(1)}\beta_j)\lambda(w w_{(4)}b_{k,(2)})\\
&&=\sum_{k,w,i,j}\lambda (S(z) w_{(1)})\lambda(w_{(3)}\alpha_i
a_kS(\alpha _j))\lambda
(z w_{(2)}\beta_i b_{k,(1)}\beta_j)\lambda(w w_{(4)}b_{k,(2)}).
\eeq
The criteria established in \ref{gencn} and
\ref{cm} imply that
$$
\sum_{k,w,i,j}\lambda (S(z) w_{(1)})\lambda
(z w_{(2)}\beta_i b_{k,(1)}\beta_j)w_{(3)}\alpha_i
a_kS(\alpha _j)\otimes w w_{(4)}b_{k,(2)}\in C^2.
$$
Hence from proposition \ref{Sinv} it follows that
\beq
&l.h.s&=\sum_{k,w,i,j}\lambda (S(z) w_{(1)})\lambda
(z w_{(2)}\beta_i b_{k,(1)}\beta_j)\lambda^S(w_{(3)}\alpha_i
a_kS(\alpha _j))\lambda^S(w w_{(4)}b_{k,(2)}).
\eeq
Now the $S$--invariance of $[z]$ together with the fact that
$\hat{\delta}([z],[z])=0$ imply that
\begin{gather*}
 l.h.s=\sum_{k,w,i,j,z}\lambda (z_{(1)} w_{(1)})\lambda
(z z_{(2)} w_{(2)}\beta_i b_{k,(1)}\beta_j)\lambda ^S(w_{(3)}\alpha_i
a_kS(\alpha _j))\lambda ^S(w w_{(4)}b_{k,(2)})\\
=
\sum_{k,w,i,j,z}\lambda (z w_{(1)})\lambda
(z \beta_i b_{k,(1)}\beta_j)\lambda ^S(w_{(2)}\alpha_i
a_kS(\alpha _j))\lambda ^S(ww_{(3)}b_{k,(2)})\hspace{.7cm}\\
=
 \sum_{k,w,i,j}\lambda
(z \beta_i b_{k,(1)}\beta_j)\lambda (z S((\alpha_i
a_kS(\alpha _j))_{(1)}))\lambda ^S(w_{(1)}(\alpha_i
a_kS(\alpha _j))_{(2)})\lambda ^S(ww_{(2)}b_{k,(2)})\\
=
\sum_{k,w,i,j}\lambda
(z \beta_i b_{k,(1)}\beta_j)\lambda ^S(z (\alpha_i
a_kS(\alpha _j))_{(1)})\lambda ^S(w_{(1)}(\alpha_i
a_kS(\alpha _j))_{(2)})\lambda ^S(ww_{(2)}b_{k,(2)})\\
=
\sum_{k,w,i,j}\lambda
(z \beta_i b_{k,(1)}\beta_j)\lambda (z (\alpha_i
a_kS(\alpha _j))_{(1)})\lambda(w_{(1)}(\alpha_i
a_kS(\alpha _j))_{(2)})\lambda(ww_{(2)}b_{k,(2)}),\hspace{0.5cm}
\end{gather*}
where  the last two equalities follow from \ref{gencn},
\ref{cm} and \ref{Sinv}.
At this point we use the fact  that $[w]\in \T$  and obtain:
\beq
&l.h.s.&= \sum_{k,i,j}\lambda
(z \beta_i b_{k,(1)}\beta_j)\lambda (z (\alpha_i
a_kS(\alpha _j))_{(1)})\lambda(w(\alpha_i
a_kS(\alpha _j))_{(2)})\lambda(w b_{k,(2)})\\
&&=
 \sum_{k}\lambda (z a_{k,(1)})\lambda(w
a_{k,(2)})\lambda (z  b_{k,(1)})\lambda(w b_{k,(2)})\\
&&=
 \sum_{k}\lambda ^S (z a_{k,(1)})\lambda ^S(w
a_{k,(2)})\lambda ^S(z  b_{k,(1)})\lambda ^S(w b_{k,(2)})\\
&&=
\sum_{k}\lambda ^S((z\star w)a_k)\lambda ^S((z\star w)b_k)=
\sum_{k}\lambda ((z\star w)a_k)\lambda ((z\star w)b_k).
\eeq
Together with the fact that $\Lambda \in \T$, this implies that $\T$ is
a monoid with respect to the $\star$--product structure as well.

It is left to show that $\T$ is invariant under the action of
$\hat{J}$, i.e
for any  $[z]\in \T$, and $\sum_ka_k\otimes b_k\in C^2$,
$$
\sum_k\lambda (J(z)a_k)\lambda (J(z)b_k)=
\sum_{k,J(z)}\lambda (J(z)_{(1)}a_k)\lambda (J(z)J(z)_{(2)}b_k).
$$
For the left hand side  one has
\beq
&l.h.s.&=
\sum_{i,j,n,m,k}\lambda(z\beta_j\alpha_i)\lambda(z\beta_m\alpha_n)\lambda(
\alpha_j\beta_ia_k)
\lambda(\alpha_m\beta_nb_k)\\
&&=
\sum_{i,j,n,m,k}\lambda(\beta_jz(\beta_mz\alpha_n)_{(2)}\alpha_i)
\lambda((\beta_mz\alpha_n)_{(1)})\lambda(\alpha_j\beta_ia_k)
\lambda(\alpha_m\beta_nb_k)\\
&&=
\sum_{i,j,n,m,k}\lambda(z z_{(2)}\alpha_{n,(2)}\alpha_i\beta_j\beta_{m,(2)})
\lambda((z\alpha_n\beta_m)_{(1)})\lambda(\beta_ia_k\alpha_j)
\lambda(\beta_nb_k\alpha_m).
\eeq
Hence, from the fact that $[z]\in \T$ and \ref{quasitriangular} (b)
and (c), it follows that
\beq
&l.h.s.&=
\sum_{i,j,n,m,k}\lambda(z(\alpha_n\beta_m)_{(1)})
\lambda(z \alpha_{n,(2)}\alpha_i\beta_j\beta_{m,(2)})
\lambda(\beta_ia_k\alpha_j)
\lambda(\beta_nb_k\alpha_m)\\
&& =
\sum_{i,j,n,m,k}\lambda(z(\beta_m\alpha_n)_{(1)})
\lambda(z \beta_j\beta_{m,(2)}\alpha_{n,(2)}\alpha_i)
\lambda(\alpha_j\beta_ia_k)
\lambda(\alpha_m\beta_nb_k)\\
&& =
\sum_{i,j,n,m,k,n',m'}\lambda(z\beta_{m'}\alpha_n)
\lambda(z \beta_j\beta_{m}\alpha_{n'}\alpha_i)
\lambda(\alpha_j\beta_ia_k)
\lambda(\alpha_m\alpha_{m'}\beta_n\beta_{n'}b_k)\\
&& =
\sum_{i,j,n',m,k}
\lambda(z \beta_j\beta_{m}\alpha_{n'}\alpha_i)
\lambda(\alpha_j\beta_ia_k)
\lambda(\alpha_m J(z)\beta_{n'}b_k)\\
&&=
\sum_{k,J(z)}\lambda (J(z)_{(1)}a_k)\lambda (J(z)J(z)_{(2)}b_k).
\eeq
\end{ssec}

\section{  K--links and K--tangles}

\begin{figure}[ht!]
\setlength{\unitlength}{1cm}
\begin{center}
\begin{picture}(6,2)
\put(0,0){\fig{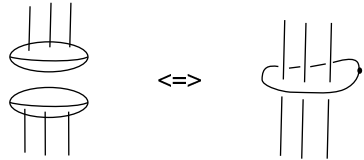}}
\end{picture}
\end{center}
\caption{Representation of 1--handle with 2--handles which pass over it}
\label{onehandle}
\end{figure}

Let $M$ be an oriented 4--dimensional manifold
together with a
  decomposition as a handlebody with a single 0--handle and a number of
  1-- and 2--handles.
Then $M$ can be represented by describing the attaching maps of the 1--
and 2--handles in $S^{3}$ \cite{K:book, 4mani}. The attaching map of
a 1--handle is
 a pair of 3--balls in $S^3$ or equivalently it can be
described as a unknot of framing 0 in $S^3$ (
figure \ref{onehandle}). In this last case the result of
attaching the 1--handle is
being thought as the manifold obtained by pushing into $B^4$ the disk
bounded by the unknot and removing a neighborhood of it. We will use
the second method putting a dot on the unknot to indicate that it
describes a 1--handle.  Then the attaching maps of the 2--handles are
described by framed links in the 1--handlebody, where if a 2--handle goes
over a
1--handle, the corresponding link component
is drown to go through the  dotted circle describing the 1--handle.

\begin{ssec} Define a Kirby link ({\em K--link})
    to be a framed link in $S^3$ where
some of the unknotted components of framing 0, bounding disjoint
Seifert surfaces, have been dotted. Then an oriented Kirby link
({\em OK--link}) is a K--link where an orientation of each link
component has been fixed. A based oriented Kirby link ({\em
BOK--link}) is an OK--link where one has fixed numbering and based
points for the undotted components and a numbering and a
set of disjoint Seifert surfaces for
the dotted components.
\end{ssec}

\begin{figure}[ht!]
\setlength{\unitlength}{1cm}
\begin{center}
\begin{picture}(6.5,4)
\put(0,0){\fig{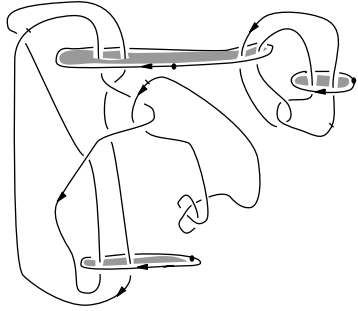}}
\end{picture}
\end{center}
\caption{A BOK--link $L$ with
 $P_L=\langle x,y,z\mid xy^{-1}xy, z^{-1}xx^{-1}z,1\rangle $}
\label{example}
\end{figure}

Given a K--link (OK--link, BOK--link)
$L$, we will denote with $M_L$ the 4--dimensional handlebody
described by $L$.
If $L$ is a BOK--link with $n$ dotted and $m$ undotted components,
then it defines a unique presentation
$\hat{P}_L=\langle x_1,x_2,\dots ,x_n \mid
 R_1,R_2,\dots ,R_m\rangle $ of $\pi_1(M_L)$, where $R_i=R_i(x_1,x_2,\dots ,x_n)$
 is a (not freely reduced) word in the $x_j$'s and shows
 in which order and with which sign the
 $i$-th undotted component intersects the Seifert surfaces of
 the dotted components starting from the base point.
   An example  is
shown in figure \ref{example}.

\begin{ssec}\label{2def}
  Two BOK--links  are said to be
  2--equivalent if and
only if they can be deformed into each other through
a sequence of the  moves (a)--(f) below (corresponding to 1-- and 2--handle
moves of the underlying 4--manifold). Changing a BOK-link through 
such a sequence will be called a 2-deformation of this link:
\begin{itemize}
\item[(a)] isotopy of framed links;
\item[(b)] any pair of one  dotted component $x$ and  one undotted
component $y$ can be removed or added if the geometric intersection
number of $y$ and the Seifert surface $S_x$ of $x$ is $\pm 1$,
while $S_x$ is disjoint from all other
dotted and undotted components (1--2 handle
cancellation or introduction);
\item[(c)] band-connected sum or difference of two
undotted link components (sliding a 2--handle over another 2--handle);
\item[(d)] band-connected sum or difference of one undotted link
component with
one dotted link component (``sliding a 2--handle over 1--handle");
\item[(e)] band-connected sum or difference of two dotted link components
(sliding an 1--handle over another 1--handle);
\item[(f)] change of numbering,  base
 points, Seifert surfaces and orientation.
\end{itemize}
The moves are illustrated in figure \ref{fig.2def}.
\begin{figure}[ht!]
\setlength{\unitlength}{1cm}
\begin{center}
\begin{picture}(13,5.5)
\put(0,0){\fig{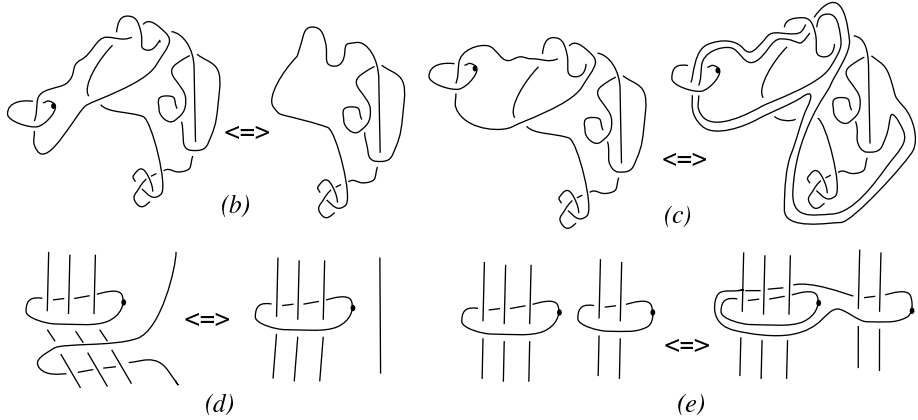}}
\end{picture}
\end{center}
\caption{Illustration of the moves (b)--(e) of a 2--deformation of
K--links}
\label{fig.2def}
\end{figure}
\end{ssec}

\begin{figure}[ht!]
\setlength{\unitlength}{1cm}
\begin{center}
\begin{picture}(10,2.5)
\put(0,0){\fig{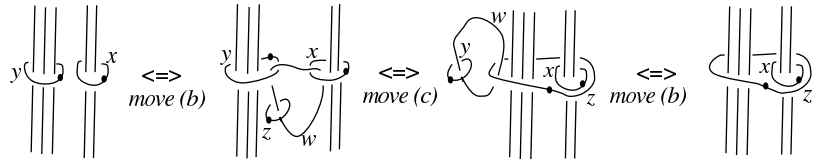}}
\end{picture}
\end{center}
\caption{Move (e) is a consequence of moves (b) and (c).}
\label{fig.prop2def}
\end{figure}

\begin{figure}[ht!]
\setlength{\unitlength}{1cm}
\begin{center}
\begin{picture}(11,2.5)
\put(0,0){\fig{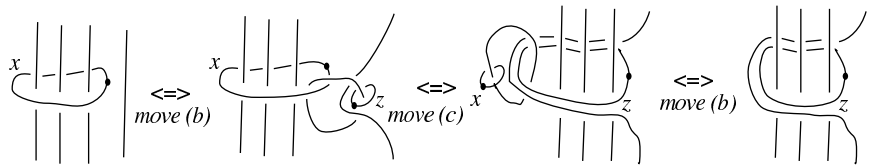}}
\end{picture}
\end{center}
\caption{Move (d) is a consequence of
    (b) and (c).}
\label{fig.prop2def1}
\end{figure}

   \begin{prop} \label{prop.2def}
    If two BOK--links can be deformed into each other
    through the moves (a)--(f) above, then they can be deformed
    into each other via moves (a), (b), (c) and (f).
    \end{prop}
   The proof is sketched in figures  \ref{fig.prop2def} and
   \ref{fig.prop2def1}.

\begin{dfn} 
Let $L$ be a BOK--link and let
$\sigma \co \hat{P}_L\rightarrow \hat{P}' $ be a sequence of AC--moves. We say that $\sigma$ can be 
lifted to $L$ if there exists a 2-deformation 
$\tilde{\sigma}\co L\rightarrow L'$ such that $\hat{P}_{L'}=\hat{P}'$.
\end{dfn}

\begin{prop}\label{connect2}
Let $L$ be a BOK--link. Then
\begin{itemize}
\item[\rm(a)] Any 2--deformation  $L\rightarrow L'$
induces a 2--deformation (sequence of AC--moves) $\hat{P}_L\rightarrow
\hat{P}_{L'}$;
\item[\rm(b)]
if
$\sigma \co \hat{P}_L\rightarrow \hat{P}'$ is a sequence of AC--moves 
then $\sigma=\xi\circ\sigma_0$, where $\sigma_0$ can be lifted to $L$ and
$\xi$ is a sequence of cancellations of terms
$x_ix_i^{-1}$ in the relations (considered as cyclic words in $x_j$'s).
\end{itemize}
\end{prop}
\begin{proof}
(a) is straightforward and
for the case when the fundamental group of the 4--thickening is
trivial, (b) is actually the statement of theorem 3.3 in \cite{SAM}.
In general one can prove (b) by induction on the length of the
sequence of AC--moves $\sigma$.
Suppose that $\sigma$ consists of a single AC--move $t$. If $t$ is not
a cancellation of a term
$x_ix_i^{-1}$ (i.e.\ the reverse direction (ii)${}^{-1}$
in  \ref{acmoves} (ii)), then it can be lifted to a single 
move $\tilde{t}\co L\rightarrow L'$,
of  type (a)$\div$(f) in \ref{2def}. 
Observe that this  is not true if $t$ is a cancellation of a term
$x_ix_i^{-1}$ in a relation, since such term implies that  the 
corresponding undotted
component enters and then goes out of the $i$-th dotted component
(without intersecting the Seifert surface of any other dotted 
component)
but possibly linking with
other undotted components or itself. Therefore, in general
we can not pull it out of the $i$-th dotted component.

 (b) will follow by induction, if we can show it for the case when 
 $\sigma=t\circ w$, where $w$ is a single AC--move of the type 
 (ii)${}^{-1}$ and $t$ is any other single AC--move (since this would 
 imply that the problematic moves can be shifted at the end of the 
sequence of AC--moves). 
Observe that if $t$ is of the type  (ii)${}^{-1}$, the statement
is trivial. 
If $t$ is of the type \ref{acmoves} (i), (iii) or (v) or (v)${}^{-1}$,
it can be easily seen that $\sigma_0$ is a single
 AC--move of the the same type as $t$ and hence we can define
$L'$ to be the BOK--link obtained by applying the move $\tilde{\sigma_0}$
on $L$.

Let $t$ be of the type
\ref{acmoves}  (iv). Suppose that the first two relations of $\hat{P}_{L}$
are  $R_1=xR_1'x^{-1}$ and $R_2$ , where
$x,y,R_1'$ are some
words in the generators, and that $w$ replaces $R_1$ with $R_1'$ 
and then $t$ replaces $R_1'$ with $R_1'R_2$.
Then define $\sigma_0$ to be the sequence of the following moves: 
conjugation of the second relation with $x$ and then multiplication of 
the first relation with the second. These moves can be lifted to $L$ 
and the resulting presentation has
as first  and second relations $xR_1'x^{-1}xR_2x^{-1}$ and
 $xR_2x^{-1}$. Obviously $R_1'R_2$, $R_2$ can be obtained from  
those  by a sequence of moves of the type  (ii)${}^{-1}$.

If $t$ is of the type
\ref{acmoves}  (ii), the only problem may arise if $R_1=xR_1'x^{-1}$,
$w$ replaces $R_1$ in $R_1'$ and then $t$ replaces $R_1'$ with $yR_1'y^{-1}$.
Then define $\sigma_0$ to be the conjugation of $R_1$ with $yx^{-1}$.
The statement follows.\end{proof}

\begin{ssec}\label{bij} We will describe 4--thickenings via their
BOK--links. In
particular, there is a surjective map $\Psi\co L\rightarrow (M_L, P_L)$
from the set of
BOK--links onto the set of  4--thickenings, 
where $\hat{P}_L\rightarrow P_L$ is
described in \ref{acmoves}. Moreover
 changing $L$ into $L'$ by 2--deformation 
moves \ref{2def} (a)$\div$ (c) and (f)  changes
$(M_L, P_L)$ into $(M_{L'}, P_{L'})$ by a 2--deformation and vice versa, i.e.\ $\Psi$ 
induces a bijection between the 2-equivalence classes of BOK--links 
onto the 2-equivalence classes of 4-thickenings. 

    Given a presentation $\hat{P}$,
with $[[\hat{P}]]$ we will denote the set of all BOK--links $L$ such
that
$\hat{P}_{L}=\hat{P}$. Suppose now that $P$ is a 2--complex realizing
$\hat{P}$ under the bijection in \ref{acmoves} and fix an element
$c\in H^2(P,Z/2)$.
Then for any $L\in [[\hat{P}]]$, $P_L=P$, and
there is an embedding
$\iota_{M_L,P}\co P\rightarrow M_L$. Denote with $[[\hat{P},c]]$ the
set of all BOK--links $L\in [[\hat{P}]]$ such that
$\iota_{M_L,P}^*(w_2(M_L))=c$. Observe that according to corollary
5.7.2 in \cite{4mani}, the
second Whitney class $w_2(M)\in H^2(M;Z/2)$ of a 4--thickening $M$,
represented by a K--link, is given by the cocycle in $H^2(M,M_1;Z/2)$
\footnote{If $M_k$ denotes the $k$-handlebody, then the boundary 
operator $H_k(M_k,M_{k-1};Z)\rightarrow H_{k-1}(M_{k-1},M_{k-2};Z)$ 
is defined by the long exact sequence on the triple 
$(M_k,M_{k-1},M_{k-2})$ and the cochain complex is obtained by dualizing 
the chain complex (see 4.2 in \cite{4mani}).}
 whose
value on each 2--handle is its framing coefficient modulo 2.
Hence,
if $\hat{P}$ has $m$ relations and $c$ is presented by a cocycle
$\bar{c}\in H^2(P,P_1;Z/2)\simeq H^2(M,M_1;Z/2)\simeq (Z/2)^{m}$,
$[[\hat{P},c]]$ is the set of all BOK--links in $[[\hat{P}]]$
whose framing coefficient on the $i$-th undotted component is equal to
$\bar{c}_i$ modulo 2.
\end{ssec}

\begin{ssec}
    We assume that the reader is familiar with the notion of
    a framed tangle, which intuitively is a slice of a framed link.
    A good reference is Shum \cite{Shum}, where it is called
    double tangle. Since
    all tangles with which we will work will be framed, 
    in the future we will just call them tangles.   A tangle with $n$
incoming and $m$ outgoing ends will be called an $n-m$ tangle.

A {\em K--tangle} will be a tangle in which some of the unknotted closed
components of framing 0, bounding disjoint Seifert surfaces, have been dotted.
An {\em OK--tangle} is a K--tangle in which an orientation of any dotted or
undotted component has been fixed, and a {\em BOK--tangle} is an OK--tangle
equipped with a choice of numbering of the closed dotted,
of the closed undotted and of
  the open components, a choice of a set of disjoint Seifert surfaces for
  all dotted components, and a choice of a
  basepoint on each undotted component $s$, where
  if the component is open, the basepoint is the positively oriented
  point in $\partial s$.

  A BOK--tangle is being described by a plane diagram which
  decomposes into a combination
of the  segments presented on
  figure \ref{fig.tangelem} and the ones obtained from them by
  changing the orientation of some components. We make the convention
  that the
  incoming ends will be drawn on the top and the outgoing ends will be
drawn on the
  bottom. The tangle plane diagrams
  used here come with a standard choice of Seifert surfaces which in the
  future won't be
  drawn, while the
  choice of base points on the closed undotted components needs to be
 indicated.

\begin{figure}[ht!]
\setlength{\unitlength}{1cm}
\begin{center}
\begin{picture}(9,4.7)
\put(0,0){\fig{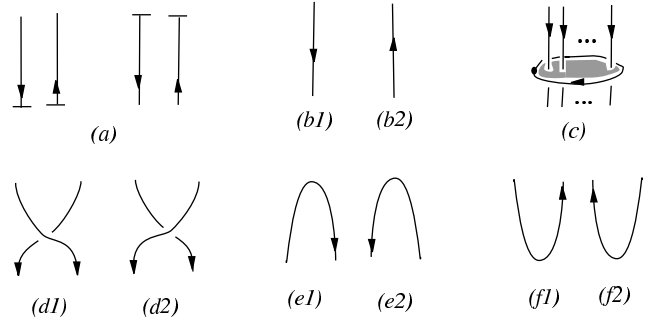}}
\end{picture}
\end{center}
\caption{Elementary tangle plane diagrams}
\label{fig.tangelem}
\end{figure}
\end{ssec}

\begin{ssec}
Two OK--tangles are equivalent
if and only if their plane diagrams can be obtained from each other via the
moves on
figure \ref{fig.RM} and \ref{fig.RM1} where  any double line
represents a number of parallel segments and
the unoriented dotted and undotted components can be oriented in any
way consistent on both sides of the identities. Two
K--tangles are equivalent if and only if their plane diagrams can be obtained
from each other via the moves on
figure \ref{fig.RM} and \ref{fig.RM1} where we have forgotten the
information about  orientation.

Observe that two K--links (i.e.\ 0-0 K--tangles) are equivalent if and only
if the corresponding framed links are isotopic.
\begin{figure}[ht!]
\setlength{\unitlength}{1cm}
\begin{center}
\begin{picture}(11,6.2)
\put(0,0){\fig{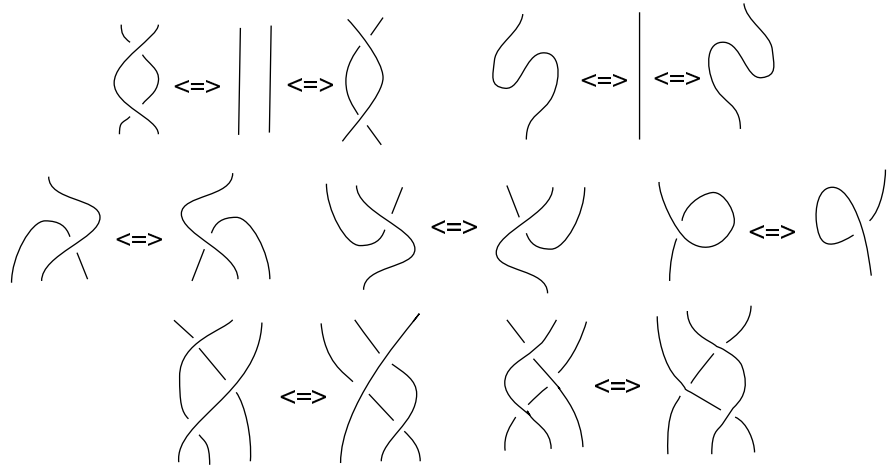}}
\end{picture}
\end{center}
\caption{``Framed'' Reidemeister moves}
\label{fig.RM}
\end{figure}

\begin{figure}[ht!]
\setlength{\unitlength}{1cm}
\begin{center}
\begin{picture}(12,9.7)
\put(0,0){\fig{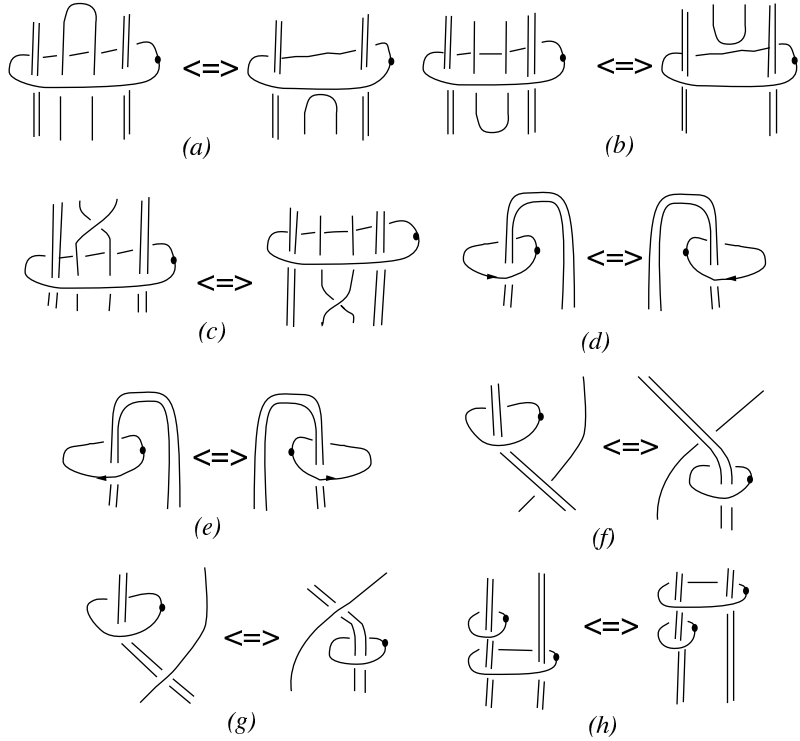}}
\end{picture}
\end{center}
\caption{Additional isotopy moves}
\label{fig.RM1}
\end{figure}
\end{ssec}

 \begin{ssec}\label{defst}
 Let $T$ be a $r-r$ K--tangle diagram with  $r$ open
  components $s_1, s_2, \ldots, s_r$  and 
  let $A_1, A_2, \ldots,A_t$ be the incoming ends and
  $B_1, B_2, \ldots,B_t$ be the outgoing
  ends of $T$ all numbered from left
 to right. Then $T$
   is called a {\em string tangle diagram} if
    there exists an element $\sigma$ in the symmetric
 group on $r$ elements ${\bf S_r}$ such that $\partial s_i=A_i\cup
 B_{\sigma (i)}$. $\sigma$ is called the underlying permutation
 of $T$. If $T$ is an OK--tangle then we add the requirement that $A_i$ 
 is the positively oriented end of $s_i$, i.e.\ the strings ``point down''.
       \end{ssec}

 \section{ Definition of the invariant}

\begin{ssec}\label{defmap}
Let $T$ be a BOK--tangle with
$n$ dotted components,
$m$ closed undotted components and $r$ open ones. Without loss of 
generality, we assume that if there are dotted components such that no 
undotted component intersects their Seifert surfaces, these are the 
first $l$ components.
By analogy with  the definition
 of the
 Hennings invariant
 \cite{Hennings}, extended to the presence of 1--handles (see for
 example in \cite{K:g}),  we define a map
  $$
  \Z(T)\co A^{\otimes (n+m)}\rightarrow A^{\otimes r}
  $$
as follows.

Let $z_j,w_i\in A$, $i=1, \ldots, n$,  $j=1, \ldots, m$.
  We refer to $z_j$ as the color of the $j$-th undotted component, and
  to $w_i$ as
  the color of the $i$-th dotted component of $T$.
  \begin{itemize}
  \item[(a)]
  Represent the BOK--tangle by  plane diagram as above;
  \item[(b)] Label the undotted components of each elementary plane diagram
  as follows:
  \begin{itemize}
      \item[-]  ``cups" and ``caps" as presented on
  figure \ref{fig.Hrules};
  \item[-]
  at each crossing of two undotted components pointing
  downwards, label
  the various segments of the  plane diagram according
  to the Hennings rules presented in figure \ref{fig.Hrules}.
  Any other crossing is obtained from those presented in the figure by
  changing the orientation of some  component $y$.  Then the
  label of $y$ changes by applying  $S^{-1}$;
   \begin{figure}[ht!]
\setlength{\unitlength}{1cm}
\begin{center}
\begin{picture}(11,2)
\put(0,0){\fig{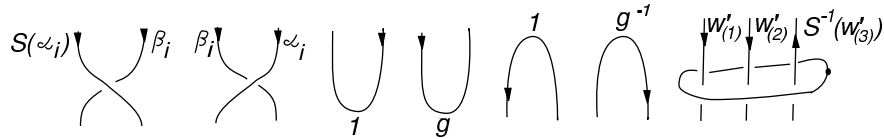}}
\end{picture}
\end{center}
\caption{Hennings type rules for labeling extended  plane diagrams}
\label{fig.Hrules}
\end{figure}

\item[-]    Let $x$ be a dotted component with color $w$ and a
Seifert surface $S_x$ and let $v_x$ be the
   normal vector of $S_x$. Let $w'=w$ if
   $v_x$ points up, and $w'=S^{-1}(w)$ if $v_{x}$ points down.
   Then, if $s_1,s_2, \ldots, s_t$ are the oriented
  segments intercepting $S_x$ , and if
  $
  \Delta ^{(t-1)}(w')=
  \sum_{w'}w'_{(1)}\otimes w'_{(2)}\otimes \ldots\otimes  w'_{(t)}
  $,
  $s_i$ gets labeled with $S^{-1}(w'_{(i)})$ if it points up,
  and with $w'_{(i)}$ otherwise as
  presented in figure \ref{fig.Hrules}.
  \end{itemize}

  \item[(c)] For each undotted component,
  starting from the base point, multiply on the right the various
  labeling elements, in the order they are
  found according to the orientation of the component. 
  In this way, one obtains an element
  $\sum_{i}a_{1,i}\otimes a_{2,i}\otimes
  \ldots a_{m,i}\otimes  b_{1,i}\otimes b_{2,i}\otimes
  \ldots b_{r,i}\in A^{\otimes (m+r)}$, where $a_{j,i}$ represents the
  product of
  the labelings of the $j$-th closed component and $b_{k,i}$
  represents the  product of
  the labelings of the $k$-th opened component. Then define
  \begin{eqnarray*}
  &&\Z(T)(z_1, \ldots,z_m,w_1, \ldots, w_n)\\
 &&=\left(\prod_{j=1}^{l}\epsilon (w_j)\right) \sum_{i} 
 \lambda ( gz_1 a_{1,i})\ldots
 \lambda ( gz_m
  a_{m,i})  b_{1,i}\otimes \ldots \otimes b_{r,i}\in A^{\otimes r}.
   \end{eqnarray*}
  \end{itemize}
  \end{ssec}

\begin{ssec}\label{observ} {\bf Remarks}\qua
    (a)\qua The application of $\epsilon
\co A\rightarrow k$ to the label of the $j$-th open component gives
exactly the invariant of the tangle $T'$ obtained from $T$
  by removing the $j$-th open component.
  \begin{itemize}
      \item[(b)] We have defined, somewhat 
     arbitrary, the value of the invariant on a disjoint 
      dotted component of color $w$ to be $\epsilon(w)$. 
      But as it   will be shown in \ref{canpr},  this is the only 
     choice consistent with  the invariance under the cancellation of 
     a dotted and undotted component (move \ref{2def} (b)). 
\end{itemize}
\end{ssec}

\begin{ssec}
We illustrate the definition with the example
of an oriented extended tangle $T$  presented
  in figure \ref{example2}.  If $w\in A$ is the
color of the dotted component and $z\in A$ is the color of the
undotted one then
$$
\Z(T)(z,w)=\sum_{i,w} \lambda(gzg^{-1}w_{(2)}\alpha _ig^{-1}w_{(1)}\beta
_i)S^{-1}(w_{(3)})\in A.
$$
 \begin{figure}[ht!]
\setlength{\unitlength}{1cm}
\begin{center}
\begin{picture}(3,2.2)
\put(0,0){\fig{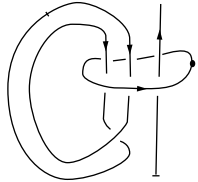}}
\end{picture}
\end{center}
\caption{An example of a BOK--tangle}
\label{example2}
\end{figure}

In the future, if we want to investigate the value of $\Z(T)$ for some
particular color of dotted or undotted component, this color may be
indicated on the plane diagram in a circle attached to the corresponding
component as in figure \ref{centr} below.
\end{ssec}

\begin{ssec} \label{proof1}{\bf Proof of theorem \ref{theo1}}\qua
The map defined so far obviously depends on the choices of numbering, base
points
and orientations. So we will start putting restrictions on the values of
  the colors in order to reduce this dependence and eventually obtain
  an invariant of BOK--links under the 2--deformation moves in
  \ref{2def}.

  The proof consists of
  showing the following  statements:
  \begin{itemize}
  \item[(A)] $\Z(T)\co Z(A)^{\otimes (n+m)}\rightarrow A^{\otimes r}$
  does not depend on the choice of base points and it is invariant
  under the moves of figures \ref{fig.RM}, \ref{fig.RM1};
  \end{itemize}
  Let now $L$ be an BOK--link with $n$ dotted and $m$ undotted
  components. Then
  \begin{itemize}
   \item[(B)]
   $\Z(L)\co Z(A)^{\otimes (n+m)}\rightarrow k$
   factors through a map $\hat{Z}(A)^{\otimes (n+m)}\rightarrow k$
   which will be denoted in the same way;
   \item[(C)]
    $\Z(L)\co \hat{Z}^S(A)^{\otimes (n+m)}\rightarrow k
   $ doesn't depend on the choice of orientation of the components
   of the link;
   \item[(D)] Let $x$ be the first, and $y$ be the second undotted
   component of $L$. Let
   also $L'$ is being obtained from $L$ by replacing $y$
   with a band connected sum of $x$ and $y$.
   Then if
   $[z],[w] \in \hat{Z}^S(A)$ are such that
     $\hat{\delta}([w],[z])=0$, and $[c]\in \hat{Z}^S(A)^{\otimes
     (n+m-2)}$, we have
   $$\Z(L)([z]\otimes [w] \otimes [c])=\Z(L')([z]\otimes [w]\otimes
   [c]).
   $$
   \item[(E)] For
   $[z],[w] \in \hat{Z}^S(A)$ let $\Z_{[z]}^{[w]}(L)$ denote the value of
   $\Z(L)$ where
   any undotted component is colored by $[z]$ and any dotted
   component is colored by $[w]$. Then  if $[zw]=[\Lambda]$,
   $\Z_{[z]}^{[w]}(L)$ is invariant under move \ref{2def} (b).
   Moreover if $[z]\in \T^4$  and $[zw]=[zw']=[\Lambda]$, then
   $\Z_{[z]}^{[w]}(L)=\Z_{[z]}^{[w']}(L)$. This common value will be
   denoted with $\Z_{[z]}(L)$.
   \end{itemize}
\end{ssec}

\begin{ssec} {\bf Proof of (A)}\qua
 First we remind Hennings' result (\cite{Hennings}) that if the colors of
the undotted
  components are in the center of the algebra,
   $\Z(T)\co  Z(A)^{\otimes
m}\otimes A^{\otimes n}\rightarrow A^{\otimes r}$ is independent of the
choice of base
points on the closed undotted components, and it is an invariant
under the  moves presented in figure \ref{fig.RM}. Moreover, from the
defining identity \ref{quasitriangular} (a) for the R--matrix
 and  the defining property of $g$, it is easy to see that
it is also an invariant  under the moves (a)$\div$ (c) on
  figure \ref{fig.RM1}.

 Suppose now that  the colors of the dotted components are in
 the center of
the algebra as well. Then the identities (f), (g) and (h) are automatically
satisfied. So, it is left to show that in this case 
(d) and (e) are satisfied as well. 
Let $x$ be the dotted component  which we want to slide
over the cup, and let $w\in Z(A)$ be its color. Since $w$ is in the
center of a ribbon algebra, $S^{2}(w)=w$ and  
(d) and (e) become equivalent. So it is enough to show (e). Let
$w'=S^{-1}(w)$. Suppose that  $n$ undotted segments pass through $x$.
Then, depending on its orientation, under the move (e) the label of the
$i$-th  segment changes as
$
g^{-1}w_{(i)}\rightarrow
S^{-1}(w'_{(n-i)})g^{-1}$ or
$S^{-1}(w_{(i)})\rightarrow w'_{(n-i)}$. But
from  \ref{propH} (d) it follows that 
$S^{-1}(w'_{(n-i)})g^{-1}=g^{-1}S(w'_{(n-i)})=g^{-1}w_{(i)}$ and
$w'_{(n-i)}=S^{-1}(w_{(i)})$. 
\end{ssec}

\begin{ssec}{\bf Proof of (B)}\qua
  The proof is based on the following observation which
  is a version of the centrality result of
  the HKR--invariant in \cite{KRS}.

      Let $T$ be a $k-l$ BOK--tangle with $n+m$ closed and
  $r$ open components. Let also
  $T'$ be the
  BOK--tangle obtained from $T$ by embracing all incoming ends
  (figure \ref{centr} (a))
  with a dotted component $x'$, and let $T''$ be  the
  BOK--tangle obtained from $T$ by embracing all outgoing ends
  with a dotted component $x''$ (figure \ref{centr} (b)).
  Fix the colors of $x'$ and $x''$ to be the same element $a\in A$ and let
  $c\in Z(A)^{\otimes (n+m)}$ describe the coloring of the
  closed components of $T$. Then
  $$
  \Z(T')(a\otimes c)=\Z(T'')(a\otimes c).
  $$

     \begin{figure}[ht!]
\setlength{\unitlength}{1cm}
\begin{center}
\begin{picture}(6,2.5)
\put(0,0){\fig{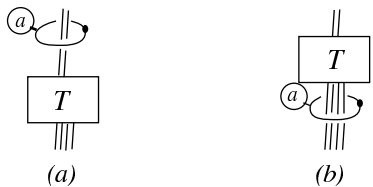}}
\end{picture}
\end{center}
\caption{Centrality of the invariant}
\label{centr}
\end{figure}

   This can be seen by decomposing the plane diagram of $T$
   into  slices such that each slice contains
   only one
   subdiagram of the type crossing, cup, cap or
   dotted component. Then, since all colors of the components of $T$
   are in $Z(A)$, one can use moves (a), (b), (c) and (h) to slide
   the dotted component colored by $a$ through.

  The statement above implies that if $T$ is an $r-r$
  string tangle then
  $
  \Z(T)$ sends $Z(A)^{\otimes (n+m)}$ into $C^r$.
  In particular, if $T$ is a 1-1 
   BOK--tangle with $(n+m)$ closed components,
   $\Z(T)$ sends $Z(A)^{\otimes (n+m)}$ into $Z(A)$.

Now we can show (B).
Let $K(A)\subset Z(A)$ be the null space of the pairing on $Z(A)$
induced by $\lambda$ as in \ref{centq}. Suppose  that
an undotted component $y$ of $L$ has a color $z\in K(A)$. Then we can use
isotopy moves to present $L$ as a closure of a 1-1 string tangle $T$
on $y$ and
$
\Z(T)$ sends $Z(A)^{\otimes (n+m-1)}$ into $Z(A)$.
Hence for any $a\in Z(A)^{\otimes (n+m-1)}$, $\Z(L)(z\otimes a)=\lambda (z\,
\Z(T)(a))=0$ by the definition of $K(A)$.

Now suppose that a dotted component $x$ of $L$ has a color $w\in
K(A)$. Since $w=w\star \Lambda$ without changing the value of the
invariant we can introduce an
undotted unknotted component $y$ of color $S(w)$  which passes once
through $x$ and in the same time change the color of $x$ to $\Lambda$
as shown in figure \ref{fig.propq1}.  But since the new tangle has an
undotted component of color $S(w)\in K(A)$ its
invariant is 0 as shown previously.
     \begin{figure}[ht!]
\setlength{\unitlength}{1cm}
\begin{center}
\begin{picture}(7,2)
\put(0,0){\fig{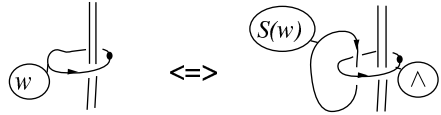}}
\end{picture}
\end{center}
\caption{Replacing a dotted component of color
$w$ with a pair of dotted component of color $\Lambda$ and
 undotted component of  color $S(w)$}
\label{fig.propq1}
\end{figure}
\end{ssec}

\begin{ssec} {\bf Proof of (C)}\label{orinv}\qua
    Observe that changing the orientation
 of a dotted
component $x$ with color $[w]\in \hat{Z}(A) $ has the
same effect as leaving its orientation the same but changing its
color to $[S(w)]$ or $[S^{-1}(w)]$. Hence if $[w]\in \hat{Z}^S(A) $, the
value of $\Z(L)$ remines unchanged.

The fact that changing the orientation of an undotted component
doesn't change the invariant is a modification of  Hennings'
argument when there is no dotted components. The link plane diagram can be
 deformed via the  regular isotopy moves of figures
\ref{fig.RM},\ref{fig.RM1} and if
 necessary changing orientation of dotted components into one
which is composed totally of segments of the types presented on
figure \ref{fig.inverse}.
     \begin{figure}[ht!]
\setlength{\unitlength}{1cm}
\begin{center}
\begin{picture}(10,4.5)
\put(0,0){\fig{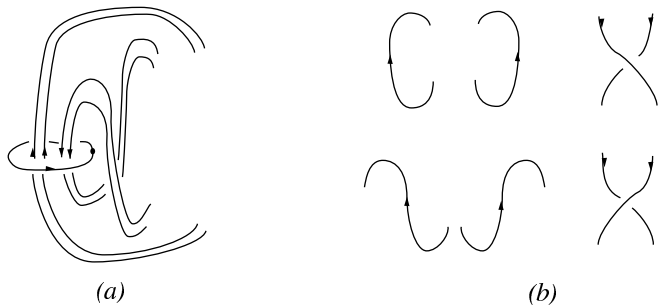}}
\end{picture}
\end{center}
\caption{Elementary plane diagrams}
\label{fig.inverse}
\end{figure}
We do this by first pulling all dotted components on the left of the
plane diagram using the moves (f) and (g) of figure \ref{fig.RM1}. In this 
way,
on the right there is left a tangle $T$ which gets closed through the
dotted components as shown in figure \ref{deform} (a).
Then, using move (c) of figure  \ref{fig.RM1} we pull all undotted
segments, which pass through a dotted component and point down, to the
right and absorb the resulting crossings into $T$ obtaining another
tangle $T'$ as shown in figure \ref{deform} (b). Then we pull down the
upper ends and pull up the lower ends of these undotted segments which
point down as they pass through a dotted component. In this way the
plane diagram is presented as the closure (through the dotted components) of
a string tangle $T''$ with positively oriented ends as shown in figure
\ref{deform} (c). At the end, by local deformations as the one on
figure \ref{deform} (d) we obtain a plane diagram in which all crossings
have the two segments pointing down.
After doing some moves  of the type of the second one in figure 
\ref{fig.RM},
we can assume that the segments of the undotted
components in $T''$ between crossings and end points are of the
type presented in figure \ref{fig.inverse} (b).
     \begin{figure}[ht!]
\setlength{\unitlength}{1cm}
\begin{center}
\begin{picture}(11,10.3)
\put(0,0){\fig{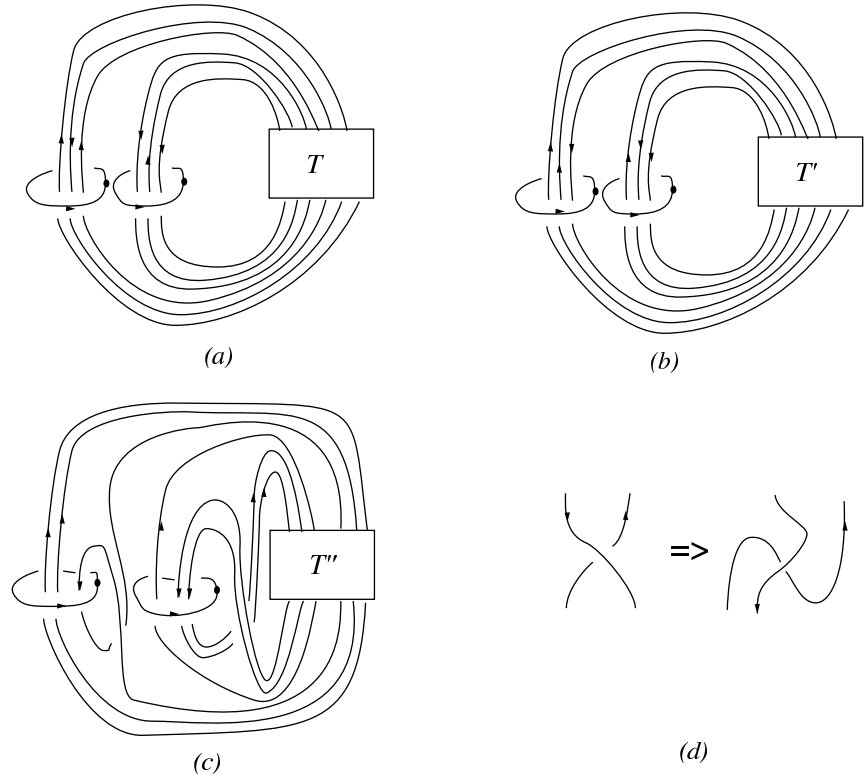}}
\end{picture}
\end{center}
\caption{Deformation of a link plane diagram}
\label{deform}
\end{figure}
Now  we want to show that, under a change of
orientation, the label of an undotted component changes by
application of $S^{-1}$.
By definition,  this is the case if we change the
orientation of an undotted component in one of the segments presented on
figures \ref{fig.inverse} (b).
Then it is enough to show the same
statement for the undotted components in figures
\ref{fig.inverse} (a). The labeling of an
undotted component which points up as it passes through a
 dotted circle of color $w$ is of the type $a=S^{-1}(w_{(i)})g^{-1}$ 
 and, after
its orientation has been changed, becomes $w_{(i)}g=gS^{-2}(w_{(i)})=
S^{-1}(a)$. The label of an
undotted component  which points down as it passes through a
 dotted circle of color $w$ is of the type $b=\alpha
 _{j,(k)}w_{(i)}S(\beta _{j,(k)})$,
and after a change of the orientation, it becomes $\beta _{j,(k)}
S^{-1}(w_{(i)})g^{-1}S(\alpha _{j,(k)})g=S^{-1}(b)$.
Since $\lambda_{gz}\circ S=\lambda_{gz}$, the statement follows.
\end{ssec}

\begin{ssec}
{\bf Proof of (D)}\qua
 First,
using isotopy moves, deform the
link plane diagram as the closure of a tangle $T$ on $y$ and $x$, where $x$
is oriented downwards and $y$ is oriented upwards as shown in figure
\ref{slide2} (a).
Without loss of
generality, we may assume that the band connected sum is like the one
presented on \ref{slide2} (b). Let
$\Z(T)([c])=\sum_i a_i\otimes
b_i\in A\otimes A$. Then,
$$
\Z (L)([z]\otimes[w]\otimes [c])=\sum_i\lambda (za_i)\lambda (wb_i).
$$
On another hand,
$
\Z(L')([z]\otimes [w]\otimes [c])=\sum_i\lambda (za_{i,(1)})\lambda (w b_i
a_{i,(2)})
$. Moreover, $a_{i,(1)}\otimes b_i a_{i,(2)}\in C^2$, since
it represents the invariant of a 2-2 string tangle. Hence,
\begin{eqnarray*}
&&\sum_{i,a_i}\lambda (za_{i,(1)})\lambda ( w b_i a_{i,(2)})=
\sum_{i,a_i,z}\lambda (z_{(1)}a_{i,(1)})\lambda (w z_{(2)}b_i a_{i,(2)})\\
&& =\sum_{i,a_i,z}
\lambda(a_{i,(1)}z_{(1)})\lambda(wb_i a_{i,(2)}z_{(2)})=
\sum_i\lambda (za_i)\lambda (w b_i ).
\end{eqnarray*}
 \begin{figure}[ht!]
\setlength{\unitlength}{1cm}
\begin{center}
\begin{picture}(6,2.5)
\put(0,0){\fig{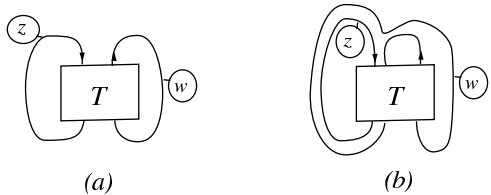}}
\end{picture}
\end{center}
\caption{On the proof of \ref{proof1} (D)}
\label{slide2}
\end{figure}
\end{ssec}

\begin{ssec} {\bf Proof of (E)}\label{canpr}\qua
    The invariance under the cancellation of a pair of dotted and undotted 
    component (move \ref{2def} (b)) is a straightforward consequence of
    the definition of $\Lambda$
and the fact that $\lambda(\Lambda)=1$ with the exception of
the case when   $L=L'\sqcup K$, where  
$K$ is a dotted 
component whose Seifert surface is disjoint from the rest of the 
link, and we have added a 
cancelling pair of dotted and undotted components such that the new
undotted component passes through $K$, obtaining in this way a new BOK--link $L''$.
Then by definition 
$\Z_{[z]}^{[w]}(L)=\epsilon(w)\Z_{[z]}^{[w]}(L')$. 
On another hand,
since $[zww]=\epsilon(w)[zw]$,  
$\Z_{[z]}^{[w]}(L'')=\epsilon(w)\Z_{[z]}^{[w]}(L')$. Hence 
$\Z_{[z]}^{[w]}(L'')=\Z_{[z]}^{[w]}(L)$ as requested.
     \begin{figure}[ht!]
\setlength{\unitlength}{1cm}
\begin{center}
\begin{picture}(10,2.7)
\put(0,0){\fig{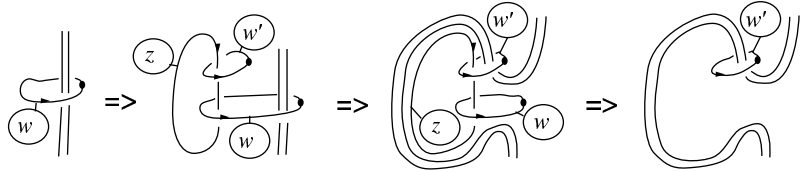}}
\end{picture}
\end{center}
\caption{On the proof of \ref{proof1} (E)}
\label{zwl}
\end{figure}

Assume now that $[z]\in \T^4$ and $[w],[w']\in \hat{Z}^S(A)$ are such that
$[zw]=[zw']=[\Lambda]$. Starting with
$\Z_{[z]}^{[w]}(L)$
we will show that one can change the color of all dotted components
from $[w]$ to $[w']$ without changing the value of the invariant.
Suppose that $x_1$ is a dotted component of color $[w]$.
Since $[w]=[\Lambda \star w]=[(z w') \star w]$, we can add a
canceling pair
of dotted component $x_2$ of color $[w']$ and an undotted component $y$ of color
$[z]$ which passes once through $x$ as shown in figure \ref{zwl}. Then,
using \ref{proof1} (D),  slide the components which pass through
$x_1$ over $y$ and since $[(z w) \star w']=[w']$,  cancel the
pair $x_1, y$. Now (E) follows from the fact that $[S(w')]=[w']$.
\end{ssec}

  We have shown that
  $\Z_{[z]} (M_L)=\Z_{[ z]} (L)$ defines an invariant of 4--thickenings. 
 To complete the proof of theorem \ref{theo1} it is left to observe 
 that $S^2\times D^2$ is represented by an undotted unknot of framing 0 
 and hence $\Z_{[z]} (S^2\times D^2)=\lambda(z)$, while $S^1\times D^3$ 
 is represented by one dotted component 
 and hence $\Z_{[z]} (S^1\times D^3)=\epsilon(w)$.
 
 Observe that if $[z]\in \T^4$, and $[zw]=[\Lambda]$, then for any 
 unit $\gamma\in k$, $[z']=[\gamma z]\in \T^4$ and $[z'w']=[\Lambda]$
 where $[w']=\frac{1}{\gamma}[w]$. Hence
 \begin{corr}\label{propo} For any unit $\gamma \in
 k$,
 $
 \Z_{[\gamma z]} (M)=\gamma ^{\chi (M)-1}\Z_{[z]} (M),
 $
 where $\chi (M)$ is the Euler characteristic of $M$.
 \end{corr}


\begin{ssec} {\bf Factorization properties of the link invariant}\qua
 Suppose that $L=L'\sqcup L''$ is a link (without dotted components),
  and $L'$ and $L''$ are sublinks of  $L$ which 
 don't have common components. Then
 let $\Z_{[z],[w]}(L'\sqcup L'')\in k$ denote the
 value of $\Z(L)$
 where all components of $L'$ have been labeled with $[z]$ and all
 components of $L''$ have been labeled with $[w]$.
 \end{ssec}

 \begin{corr}\label{fact} \begin{itemize}
     \item[\rm(a)] If $[z],[w]\in \hat{Z}^S(A)$ are such that
 $\hat{\delta}([w],[z])=0$,
 then
     $\Z_{[w],[J(z)]}(L'\sqcup L'')=\Z_{[w]}(L')\Z_{[J(z)]}(L'');$
     \item[\rm(b)] If $[z]\in \T$ then 
      $\Z_{[z\star J(z)]}(L)=
     \Z_{[z]}(L)\Z_{[J(z)]}(L)$.
     \end{itemize}
     \end{corr}
      \begin{figure}[ht!]
\setlength{\unitlength}{1cm}
\begin{center}
\begin{picture}(7,1.5)
\put(0,0){\fig{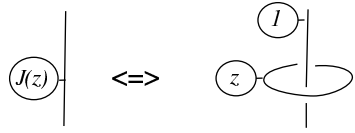}}
\end{picture}
\end{center}
\caption{Replacing an undotted component of color
$J(z)$ with an undotted component of color $1$ embraced by an undotted
component of color $z$}
\label{Jcol}
\end{figure}
\begin{proof}
    The definition of $J$ in \ref{Jdef} implies that coloring a
     component $x\in L''$ with $[J(z)]$ is equivalent to coloring $x$
     with 1 and embracing it
     with a small undotted unknot $x'$ of color $[z]$
     as showed in figure \ref{Jcol}.
     But since any component $y\in L'$ has color $[w]$,
     according to \ref{proof1} (D),
      $y$ can be slided over $x'$ and   it is a basic fact from the
     Kirby calculus, that in this way $y$ can be unlinked from $x$.
     Hence we can unlink any component of $L'$ from any component of  $L''$ and
     move them apart. This shows (a).

     Let now $L^{\#}=L\sqcup L'$ be the double of $L$, i.e.\ $L'$ is a copy
     of $L$, and $L^{\#}$ is obtained from $L$
     by adding a parallel  to each component of $L$, using the
     framing. Then (b) would follow from (a) if we could
     show that for any $[z],[w]\in
     \hat{Z}^S(A)$,
     $$
     \Z_{[z\star w]}(L)=\Z_{[z],[w]}(L'\sqcup L).
     $$
     Let $x$ be a component of $L$ colored by $[z\star w]$.
     $L$ can be presented as a closure of a 1-1 string tangle $T$ on $x$
     with $\Z_{[z\star w]}(T)=c\in Z(A)$. Then
     $$
     \Z_{[z\star w]}(L)=\lambda ((z\star w)c)=\lambda (w(S(z)\star
     c))=\lambda(zc_{(1)})\lambda(wc_{(2)}),
     $$
     where in the last two equalities we have used \ref{strco} and \ref{Sinv}.
     But the last expression is exactly the invariant of a link
     obtained from $L$ by adding a parallel component $x'$ of $x$
     and coloring $x$ by $[w]$ and $x'$ by $[z]$.\end{proof}

\begin{ssec}\label{proof2}{\bf Proof of lemma \ref{lemma2}}\qua
Let $P$ be a 2--dimensional CW complex, $c\in H^2(P,Z/2)$, and  let
$$\hat{P}=\langle x_1,x_2,\dots ,x_n \mid R_1,R_2,\dots ,R_m\rangle .
$$
From   \ref{bij} it follows that in order to prove lemma
\ref{lemma2} it is enough to show that, if  $L_0$ is a standard
representative in  $[[\hat{P},c]]$, then for  any
other $L\in [[\hat{P},c]]$ and any $[z]\in \T^2$,
$\Z_{[z]}(L)=\Z_{[z]}(L_0)$. So, we proceed with the description of
$L_0$.

Without loss of generality we assume that, if $\hat{P}$ contains
trivial relations, these are the last $k$ relations. Then let
$$Q=R_1R_2\dots
R_{m-k}=x_{i_1}^{e_1}x_{i_2}^{e_2}\ldots x_{i_t}^{e_t}, \;\mbox{where}\;
e_i=\pm 1,
$$ be the unreduced word
obtained by putting together all nontrivial relations in $\hat{P}$.
Let also $t_i^{+}\;( t_i^{-})$ denote 
the absolute value of the sum of the
 positive (negative) exponents of $x_i$ in $Q$ and $l_i$ denote the
  length of
the relation $R_i$  (the sum of the absolute
values of the exponents of $x_j$'s in $R_i$).
Define  $\sigma _Q $ to
 be the permutation element in the symmetric group  ${\bf S_t}$, such that
 $
 \sigma_Q(k)< \sigma_Q(l)$ if ($i_k< i_l$)  or 
 ($i_k=i_l$  and $e_k< e_l$) or ($i_k=i_l$, $e_k=e_l$ and $k<l$).
 Observe that applying the permutation $\sigma_Q$ on the letters of $Q$
 gives the word
 $x_1^{-t_1^-}x_1^{+t_1^+}x_2^{-t_2^-}x_2^{t_2^+} \ldots
x_n^{-t_n^-}x_n^{t_n^+}$.
Let also $\tau_Q$ be the following element in ${\bf S_t}$ presented as
product of cycles:
\beq
&\tau_Q=&(\sigma_Q(1),\sigma_Q(2), \ldots,\sigma_Q(l_1))
  (\sigma_Q(l_1+1),\ldots,\sigma_Q(l_1+l_2))\ldots\\
&&\ldots(\sigma_Q(t-l_m+1),\ldots, \sigma_Q(t)).
\eeq
Fix a braid $B_Q$ on $t$ strings oriented downwards, which has $\tau
_Q$ as  underlying permutation.
Then the standard representative $L_0$  is defined to be
the BOK--link in $[[P,c]]$ of the type presented in figure \ref{deform} (c),
where the dotted components are ordered in increasing order from the
right to the left and
where $T''=T_0$ is a string tangle which is obtained
by putting next to $B_Q$
 $k$ undotted unknots.  The framing coefficients of all undotted
components are chosen to be 0 or 1 depending on the corresponding value of the
cocycle $\bar{c}\in H^2(P,P_1; Z/2)$.

Let $L$ be another BOK--link in $[[\hat{P},c]]$
and let $[z]\in \T^2$, i.e.\ there exist $[z_1],[z_2]\in \hat{Z}^S(A)$
such that $[z]=[z_1J(z_2)]$ and $\hat{\delta}([z_1],[z_2])=0$. 
By the definition of $\hat{P}_L$,
each letter $x_{i_j}^{e_j}$ in $Q$ corresponds to an
intersection point $A_j$ in $L$
of an undotted component with the Seifert surface of
the $i_j$-th dotted component and $e_j$ is the sign of this
intersection. Use $\sigma_Q$ to define an order of the set of points 
$A_j$, in particular
$A_j\prec A_l$ if $\sigma_Q(k)< \sigma_Q(l)$.
By isotopy moves as  in \ref{orinv},
we  deform $L$ into a link $L''$ from the type
presented in figure \ref{deform} (c) so that
the points $A_j$ are ordered in increasing order from the
right to the left. Then
$\Z_{[z]}(L)=\Z_{[z]}(L'')$. Of course, $T''$ in general will be different
from $T_0$, but it is a string tangle and since $\hat{P}_L=\hat{P}_{L_0}$, 
 $T''$ has the same underlying
permutation $\tau_Q$.
Now, as shown in figure \ref{Jcol}, labeling an undotted
component $y\in L''$
with $[z_1J(z_2)]$ is the same as labeling $y$ with $[z_1]$ and embracing
it with undotted component $y'$ labeled by $[z_2]$. But since 
$\hat{\delta}([z_1],[z_2])=0$, any component
labeled by $[z_1]$ can be slided over any component labeled by $[z_2]$.
Therefore if $x$ is any other undotted component in $L''$,
we can use sliding of $x$ over $y'$ to
change the sign of any crossing of $y$ with $x$ and
 by sliding $y$ over $y'$ we can
add two positive or two negative twists
on $y$, i.e.\ change the framing coefficient of $y$ with $\pm 2$.
 Since $T''$ and $T_0$
have the same underlying permutation, by applying a
sequence of such operations $T''$ can be transformed into $T_0$.
Hence $\Z_{[z]}(L_0)=\Z_{[z]}(L'')
=\Z_{[z]}(L)$.
\end{ssec}

\begin{ssec}{\bf Proof of corollary \ref{corr2}}\qua
  If $A$ is a
finite-dimensional unimodular triangular Hopf algebra
then  the positive and negative crossings of two undotted components 
have the same labeling. Therefore for any 
$[z]\in\T^4$ we can repeat the argument above and show that if 
$L_1,L_2\in [[\hat{P},c]]$ then $\Z_{[z]}(L_1)=\Z_{[z]}(L_2)$.
 Moreover, the ribbon element in a triangular algebra is $\theta=1$.
 Hence the invariant in
lemma \ref{lemma2} won't depend any more on the framings of the
undotted components, in particular for any $L_1,L_2\in  [[\hat{P}]]$ and
 any $[z]\in \T^4$,
 $$
 \Z_{[z]}(L_1)=\Z_{[z]}(L_2).
$$
Now, let
$\hat{P}\rightarrow \hat{P'}$ be an AC--move and $L\in [[\hat{P}]]$.
By \ref{connect2} (b)  there exists  a BOK--link $L'$,
 2--equivalent to $L$ such that $\hat{P}'$
can be obtained from  $\hat{P}_{L'}$ by cancellations of terms of
the type $x_{i}x_i^{-1}$.
But such term in $L'$ corresponds to an undotted segment which enters into
the $i$-th dotted component $x_i$, possibly links with other undotted components 
or itself (but
doesn't pass through other dotted ones) and then goes out of $x_i$. Now
by cross changes we can unlink any such undotted component and then by
isotopy moves, pull it
out of $x_i$ without changing the value of the invariant.
The result is an BOK--link $L''\in [[\hat{P}']]$ and we have
$\Z_{[z]}(L'')=\Z_{[z]}(L')=\Z_{[z]}(L)$.
\end{ssec}

\section{Relation with the 3--manifold invariants}

\begin{ssec}\label{moves3}
Suppose that  we want an invariant of a 4--thickening
to depend only on its boundary.
This would imply (see \cite{K:book}) invariance under
two additional
moves:
\begin{itemize}
\item[(i)]  Removing or adding  a dot on an 0--framed unknot.
This corresponds to replacing a one handle with its canceling
2--handle and vice versa;
\item[(ii)] Deleting or adding an unknot $U^{\pm 1}$ of framing $\pm 1$,
contained in a neighborhood disjoint from the rest of the link, which
corresponds
to taking a connected union with $CP^2$ or $\overline{CP^2}$.
\end{itemize}
\end{ssec}

 In general, $\Z_{[ z]} $ won't be invariant  under
 these additional moves, but in many
 examples (including all the ones coming from the quantum $sl(2)$)
 $\Z_{[z]} $ can be normalized to depend
only on the boundary.  We will use the statement below only for $[z]\in
\T$, but observe that it is true in the following weaker form:

\begin{prop} \label{ccx}
    Suppose that $[z]\in \T_{Z}$
    and that $[z J(z)]=X[\Lambda]$
for some unit $X\in k$.
Then
$
X=\lambda(z\theta ^{-1})\lambda(z\theta).
$
\end{prop}
\begin{proof} Since $\epsilon(\theta)=\sum_i\epsilon(\beta_ig\alpha_i)=1$,
$[\theta^{-1}zJ(z)]=X[\Lambda]$ and therefore $X=\lambda (\theta^{-1}zJ(z))$.
Substituting here the expression for $J(z)$ from \ref{exJ} we obtain that
$$
X=\lambda (\theta^{-1}zJ(z))=
\lambda (z((S(z)\theta)\star \theta^{-1})).
$$
Now since $[S(z)]=[z]$, applying  \ref{strco}  it follows that
$$
X=
\lambda((z\theta)(z\star \theta^{-1}))=\lambda((z\theta)(1\star
z\theta^{-1}))=\lambda(z\theta^{-1})\lambda(z\theta),
$$
where in the second equality we have used the fact that 
$[z]\in\T_Z$.\end{proof}

\begin{ssec} {\bf Proof of corollary \ref{3inv}}\qua
    Since $C_{\pm }=\lambda(z\theta^{\pm 1})$, the first assertion
    follows from the proposition above.
The rest follows from the observation
that the ordered pair $(\sigma_{+}-n,\sigma_{-}-n)$ is an invariant
under  2--deformations of $M$ since a 2--handle slide \ref{2def} (c)
doesn't change the number of dotted components
and the values of $\sigma_{+}$ and  $\sigma_{-}$,
 while move \ref{2def} (b) reduces by one the
number of dotted components, and in the same time reduces by one the
values of $\sigma_{+}$ and  $\sigma_{-}$. Moreover, the
proposition \ref{ccx}
implies that under the moves \ref{moves3}(i) and (ii), $\Z_{[z]}(M)$ changes
exactly as $C_+^{\sigma_{+}-n}C_-^{\sigma_{-}-n}$ and therefore
their quotient $\Z_{[ z]}^{\partial}(\partial M)$ depends
only on the boundary.
\end{ssec}

\begin{prop} \label{3fact} Let $[z]\in \T ^3$. Then
\begin{itemize}
\item[\rm(a)] for any unit $\gamma\in k$, $[\gamma z]\in \T ^3$ and
$\Z_{[\gamma z]}^{\partial}(\partial M)=
\gamma ^{\sigma _0}\Z_{[ z]}^{\partial}(\partial M)$;
\item[\rm(b)] if $[J(z)], [z\star J(z)]\in \T ^3$ then
$\Z_{[ z\star J(z)]}^{\partial}(\partial M)=\Z_{[ z]}^{\partial}(\partial
M)\Z_{[J( z)]}^{\partial}(\partial M).$
\end{itemize}
\end{prop}

The proposition is a direct consequence of the corollaries \ref{propo}
and \ref{fact}.

\begin{corr} \label{3factfact} If $A$ is $\Lambda$--factorizable then 
for any $[z]\in\T^3$, 
$$
\Z_{[z]}^{\partial}(\partial
M)\Z_{[J(z)]}^{\partial}(\partial M)=
X_z^{\sigma_0}\Z_{[1]}^{\partial}(\partial M).
$$
\end{corr}
\begin{proof}
Since the algebra
is $\Lambda$--factorizable,  $J(1)=\gamma \Lambda $. Then 
$[z]\in\T^3$ implies that $[zJ(z)]=X_z[\Lambda]$. 
Applying 
$\frac{1}{\gamma}J$ on both sides of the equality and using \ref{corrJ}, 
we obtain that 
$[J(z)\star J^2(z)]=\gamma X_z \,[1]$. But \ref{J} (b) implies that 
$[J^2(z)]=\gamma [S(z)]=\gamma [z]$. Hence $[J(z)\star z]=\gamma X_z \,[1]$.
Since in this case $J$ is a bijection, we can reverse the argument and 
therefore obtain that,  if the algebra
is $\Lambda$--factorizable,
$$
 \T ^3=\{[z]\in \T\mid  [z\star J(z)]=X_z[1]  \mbox{ for some unit}\; X_z\in k\}.
$$
In particular, if
$[z]\in \T ^3$ then $[J(z)], [z\star
J(z)]=X_z[1]\in \T ^3$. Now the statement follows from
 proposition \ref{3fact}.\end{proof}
 

\section{Examples}

To illustrate the generality of the present framework we describe two
examples. The first one is useful to get familiar with the framework,
and the second one is the quantum $sl(2)$ case, which shows quite rich
algebraic structure, but it is not interesting for the AC--conjecture.
Indeed all
$sl(2)$ theories are actually 3--dimensional.

\begin{ssec} {\bf The cocommutative case: $R=1\otimes 1$}

Since this is a particular case of a triangular structure on $A$,
we are talking about invariants of 2--complexes. First, observe that
in this case $g=1$ and $S^2=1$. As a consequence,
the invariant has very simple definition, which is worth
writing down. Let $[z]\in \T^4$
and choose $[w]\in \hat{Z}^S(A)$ such that $[zw]=[\Lambda]$.
Let $\hat{P}=\langle x_1,x_2,\dots ,x_n \mid R_1,R_2,\dots ,R_m\rangle $ be a
presentation, where $R_i=R_i(x_1,x_2,\ldots,x_n)$. Let also $Q, 
\,\sigma_Q,\,
t_i^{\pm}, \,l_j$
and $t$ be as in \ref{proof2} and $t_i=t_i^++t_i^-$
  be the total exponent of $x_i$.
 Associated to $Q$, define a bijective map
$
S_Q\co A^{\otimes t}\rightarrow A^{\otimes t}$ such that
$$
S_Q(\sum_i a_{1,i}\otimes  a_{2,i}\otimes \ldots  \otimes a_{t,i})=
\sum_i S^{\epsilon_1}(a_{1,i})\otimes  S^{\epsilon_1}( a_{2,i})\otimes \ldots
\otimes S^{\epsilon_1}( a_{t,i}),
$$
where $\epsilon_j=(1-e_j)/2$ and $S^0=id_A$, i.e.\ in case that the $j$-th exponent
in $Q$ is negative $S_Q$ applies the antipode on the $j$-th factor in
$  A^{\otimes t}$.

Let
$\bar{\sigma}_Q\co A^{\otimes t}\rightarrow A^{\otimes t}$ be  the
permutation of factors induced by $\sigma_Q$ and let
$$
\sum_i a_{1,i}\otimes  a_{2,i}\otimes \ldots  \otimes
a_{t,i}=S_Q\circ\bar{\sigma}_Q^{-1}
(
\Delta^{t_1-1}w\otimes  \Delta^{t_2-1}w \otimes \ldots
\otimes \Delta^{t_n-1}w)\in A^{\otimes t}.
$$
Then from the definition of $\Z _{[z]}^2$ in section 7 and the
 fact that we are in the case when $R=1\otimes 1$, it follows that
\beq
&\Z_{[ z]}^2 (P)=&
\sum _{i}\lambda (z a _{1,i} a _{2,i} \ldots a_{l_1,i})
\lambda (z a_{l_1+1,i} a_{l_1+2,i} \ldots a_{l_1+l_2,i}) \ldots\\
&&\ldots
\lambda (z a_{t-l_m+1,i} a_{t-l_m+2,i} \ldots a_{t,i}).
\eeq

We illustrate the technique with the case of a group
 algebra and $[z]=1$. The result is a well known invariant which depends
 on the fundamental group of $P$.
 Let $A=k[G]$, where $G$ is a finite group. Then the product on $A$
 is induced from the one in $G$, and for any $a\in G$,
 $\Delta (a)=a\otimes a$ and
  $S(a)=a^{-1}$.
 $A$ is a unimodular algebra with $\Lambda =\sum _{a\in G}a$, and
 $\lambda\in A^*$ defined as
 $\lambda (1) =1$, and $\lambda (a)=0$ if $a\neq 1$.
 Hence the algebra is cosemisimple, and it is semisimple if and
 only if the characteristic of $k$ doesn't divide the order of $G$.
  For $z=1$ and
 $w=\Lambda$, the value of the
 invariant is:
 $$
 \Z_{[1]}^2(P)=\sum _{\{a_j\}_{j=1}^n}\lambda (
R_1(a_1,\dots ,a_n))
\lambda ( R_2(a_1,\dots ,a_n)) 
\ldots \lambda ( R_m(a_1,\dots ,a_n)),
$$
where the sum is over all possible sequences ${\{a_j\}_{j=1}^n}$ of 
elements in $G$ and
 $R_i(a_1,a_2,\dots ,a_n)$ denotes the image of the word $R_i$
under the group homomorphism of the free group on the generators
$x_1,\,x_2, \ldots,x_n$ into $G$ given by $x_j\rightarrow a_j$.
Hence $\Z_{[1]}^2(P)$ is equal to the number of all possible group
homomorphisms $G\rightarrow \pi _1(P)$.
\end{ssec}

\begin{ssec}{\bf  The quantum enveloping algebra of $sl(2)$}
\end{ssec}

  We  use here the definition of the 
finite-dimensional quantum enveloping algebra of $sl(2)$  ``at root 
of unity'' as given in chapter 36 of the book of 
G. Lusztig \cite{L:book},
and we refer the reader to  \cite{L:book},
chapters 23, 31, 32, 34 and 36, for  the proof that the definition  is 
consistent with the Hopf algebra axioms and that the category of 
representations of the algebra is the same as the one of the 
finite-dimensional quantum $sl(2)$, defined in a 
more familiar ways.
For the $sl(2)$ case, many statements can actually be
easily verified by direct computation as well.

\begin{ssec} \label{defn}
Let $p>3$ be a prime number and let $k'=Z[v]/\langle 1+v+\ldots+v^{p-1}\rangle $ and
$k=Q[v]/\langle 1+v+\ldots+v^{p-1}\rangle $.  For any $n,m\in Z$ such that
$m\geq 0$ we will use the following
common notations:
\beq
&&[n]=\frac{v^n-v^{-n}}{v-v^{-1}},\quad
\qbi{n}{m}=\frac{\prod _{s=0}^{m-1}(v^{n-s}-v^{-n+s})}{\prod
_{s=1}^{m}(v^{s}-v^{-s})},\\
&&\{m\}=\prod_{i=1}^m
(v^i-v^{-i}) ,\;\; \{0\}=1,
\eeq
hoping that the double use of square bracket to denote equivalence
classes in $\hat{Z}(A)$ and quantum integers will not bring to a
confusion. Note that $\{p-1\}=p$.
Define $A$ to be the
$k$ algebra generated by the elements $1_c E^{(n)}$, $1_cF^{(n)}$ such that
$c\in Z/p$ and $0\leq n \leq p-1$ and relations:
\begin{eqnarray*}
&&1_cE^{(n)}1_sE^{(m)}=\delta _{c,s+2n}\qbi{n+m}{n}1_cE^{(n+m)};\\
&&1_cF^{(n)}1_sF^{(m)}=\delta _{c,s-2n}\qbi{n+m}{n} 1_cF^{(n+m)};\\
&&1_cF^{(n)}1_sE^{(m)}=\delta_{c,s-2n}\sum_{t=0}^{min(m,n)}\qbi{m+n-s}{t}
1_cE^{(m-t)}1_{c-2(m-t)}F^{(n-t)};\\
&&1_cE^{(n)}1_sF^{(m)}=\delta_{c,s+2n}\sum_{t=0}^{min(m,n)}\qbi{m+n+s}{t}
1_cF^{(m-t)}1_{c+2(m-t)}E^{(n-t)}.
\end{eqnarray*}
We introduce the notation
$1_cE^{(n)}F^{(m)}=1_cE^{(n)}1_{c-2n}F^{(m)}$. Then
$A$ is a finite-dimensional  algebra with identity
$\1=\sum_{c\in Z/p}1_c$ and basis $\{ 1_cE^{(n)}F^{(m)}\}$, where $c\in
Z/p$, $0\leq
n,m\leq p-1$. $A$ has a Hopf algebra structure with the
following structure maps:
\begin{eqnarray*}
&&\epsilon(1_cE^{(n)})=\epsilon(1_cF^{(n)})=\delta _{c,0}\delta
_{n,0};\\
&&\Delta (1_cE^{(n)})=\sum_{a=0}^n\sum_{r\in
Z/p}v^{a(a-n)+r(n-a)}1_rE^{(a)}\otimes
1_{c-r}E^{(n-a)};\\
&&\Delta (1_cF^{(n)})=\sum_{a=0}^n\sum_{r\in
Z/p}v^{a(a-n)-(c-r)a}1_rF^{(a)}\otimes
1_{c-r}F^{(n-a)};\\
&&S(1_cE^{(n)})=(-1)^{n}v^{n(c-1-n)}1_{-c+2n}E^{(n)};\\
&&S(1_cF^{(n)})=(-1)^{n}v^{-n(c-1+n)}1_{-c-2n}F^{(n)};
\end{eqnarray*}
It is easy to check that $A$ is a unimodular Hopf algebra
with an integral $\Lambda =1_0E^{(p-1)}F^{(p-1)}$ and that $A^*$ has
as a right integral $\lambda$ defined as
$$
\lambda(1_cE^{(n)}F^{(m)} )=v^{c}\delta _{n,p-1}\delta _{m,p-1}.
$$
Obviously, $\lambda(\Lambda )=1$. $A$ is a quasitriangular ribbon algebra
with
$$
R=\sum _{n=0}^{p-1}\sum_{r,s\in
Z/p}v^{\frac{n(n-1)}{2}+\frac{rs}{2}}\{n\}1_{r}F^{(n)}\otimes
1_{s}E^{(n)}\;\mbox{ and }\;
g=\sum_{c\in Z/p}v^{-c}1_c.
$$
\end{ssec}

\begin{ssec}\label{eval}The center of  $A$ is  described in
\cite{K:mcg}, where the following notations are used:
$K=\sum_{s\in Z/p}v^s 1_s$, $\pi_s(K)=1_{-2s}$,
$E=(v-v^{-1})\sum_{c\in Z/p}1_cE^{(1)}$ and $F=\sum_{c\in Z/p}1_cF^{(1)}$.
Following \cite{K:mcg} we define
\beq
&&X=(v-v^{-1})\sum_{s=0}^{p-1}1_sE^{(1)}F^{(1)}+\sum_{k=1}^{p-1} b(k-1)
1_{2k}\in Z(A)\quad\mbox{and }\\
&&\phi_j(x)=\prod_{0\leq s\leq p-1:b(s)\neq b(j)}(x-b(s))\in k[x],\;
j=0,\ldots,q
\eeq
where $b(s)=b(p-1-s)=\frac{v^{2s+1}+v^{-2s-1}}{v-v^{-1}}$. Let
$q=\frac{p-1}{2}$
and let
\beq
&&P_j=\frac{1}{\phi_j(b(j))}\phi_j(X)-\frac{\phi_j'(b(j))}{\phi_j(b(j))^2}
\phi_j(X)(X-b(j)),\; j=0, \ldots, q,\\
&&N_j=\frac{1}{\phi_j(b(j))}\phi_j(X)(X-b(j)), \; j=0, \ldots, q-1, \\
&&N_j^+=T_j N_j,\; N_j^-=(1-T_j) N_j,\; \mbox{where }\;
T_j=\sum_{s=j+1}^{p-1-j} 1_{-2s}.
\eeq
Lemma 18 in \cite{K:mcg} allows to express the elements above  in terms of the
algebra basic elements $1_sE^{(i)}F^{(j)}$ as follows:
$$
1_{-2s}\phi_k(X)(X-b(k))=\sum_{j=0}^{p-1}
\prod_{i=j+1}^{p-1}(b(k)-b(i+s))([j]!)^2(v-v^{-1})^j 
1_{-2s}E^{(j)}F^{(j)},$$
$$
1_{-2s}\phi_k(X)=\sum_{j=0}^{p-2}\sum_{t=j+1}^{p-1}\prod_{i=j+1,
i\neq t}^{p-1}(b(k)-b(i+s))([j]!)^2(v-v^{-1})^j 1_{-2s}E^{(j)}F^{(j)},$$
\beq
&&\Phi_k(b(k))=([p-1]!)^2\frac{(v-v^{-1})^{p-2}}{[2k+1]^2},\\
&&\Phi_k'(b(k))=([p-1]!)^2\frac{(v-v^{-1})^{p-3}[2(2k+1)]}{[2k+1]^5},
\eeq
for any $k=0, \ldots,q-1$, and 
\beq
&&1_{-2s}\phi_q(X)=\sum_{j=0}^{p-1}
\prod_{i=j+1}^{p-1}(b(q)-b(i+s))([j]!)^2(v-v^{-1})^j1_{-2s}E^{(j)}F^{(j)},\\
&&\Phi_q(b(q))=([p-1]!)^2(v-v^{-1})^{p-1}.
\eeq
From here one can see that $N_0^-=(v-v^{-1})\Lambda$ and
$\lambda(N_i^-)=(v-v^{-1})[2i+1]^3$. In particular
$\lambda(N_i^-)\neq 0$
for any $i=0,\ldots, q-1$.
\end{ssec}

\begin{ssec}\label{th}({\bf Kerler \cite{K:mcg}})\qua
$Z(A)$ is a $3q+1$ dimensional
algebra with
basis $\{P_i,N_j^{\pm}, i=0,\ldots,q, \,j=0, \ldots,q-1\}$ and
products:
\beq
&& P_iP_j=\delta_{i,j}P_j\\
&&P_iN_j^{\pm}=\delta_{i,j}N_j^{\pm}\\
&&N_l^{\pm}N_j^{\pm}=N_l^{\mp}N_j^{\pm}=0.
\eeq
Moreover, the ribbon element in this basis is given by
$$
\theta
=v^qP_q+\sum_{j=0}^{q-1}v^{2j(j+1)}
  (P_j+\frac{2j+1}{[2j+1]}N_j-\frac{p}{[2j+1]}N_j^-).
  $$
 Observe that since $X$ and $T_j$ are $S$--invariant, any element in $Z(A)$ is
$S$--invariant and
$$
K(A)=span\{P_q, N_j, \, j=0, \ldots,q-1\}.
$$
Hence $\hat{Z}(A)=\hat{Z}^S(A)$ is generated by $[P_i],[N_j^{-}],
 \,i,j=0, \ldots,q-1$
and the following relations:
\begin{itemize}
\item[(a)] $ [P_i][P_j]=\delta_{i,j}[P_j]$,
\item[(b)] $[P_i][N_j^{-}]=\delta_{i,j}[N_j^{-}]$,
\item[(c)] $[N_l^{-}][N_j^{-}]=0$.
\end{itemize}
To be able to continue we need to understand also the $\star$ product
structure of the algebra. An easy calculation shows that
\begin{itemize}
\item[(d)] $ \hat{J}([1])=\gamma_p[\Lambda] \;\mbox{ and }\;
\hat{J}([\Lambda])=[1]=\sum_{i=1}^{q-1}[P_i]$,
\end{itemize}
where $\gamma_p=p^3$, i.e.\ the algebra is
$\Lambda$--factorizable. Then according to corollary \ref{corrJ},
$J^2=\gamma_p \,
\1$, $\gamma _p^{-1}J\co Z(A)\rightarrow Z_{\star}(A)$
is an algebra isomorphism and therefore the $\star$ algebra structure can
be derived from the knowledge of $J$.
\end{ssec}

\begin{lemma} \label{Jact} $\hat{J}([N_i^-])=(v-v^{-1})[2i+1]^2\sum
_{k=0}^{q-1}\frac{[(2i+1)(2k+1)]}{[2k+1]}[P_k]$.
\end{lemma}

We will need  the following proposition:

\begin{prop} \label{poly} For any $b$ such that $0\leq b\leq p-2$,
    let $\Omega _{b}=
Z(A)\cap span\{1_{s}E^{(a)}F^{(a)},\, s\in Z/p,\, 0\leq a\leq
b\}$.
Then $\Omega _{b}\subset span\{P_i,\,N_j,\, 0\leq i\leq q,\,
0\leq j\leq q-1\}$.
\end{prop}
\begin{proof}
We will show that $\Omega _{b}=span\{X^a,\, 0\leq a\leq b\}$.
Then the statement will follow from the observation in \cite{K:mcg} 
that any polynomial in $X$ is
contained in the span of $P_i$, $i=1,\ldots,q$ and $N_j$,
$j=1,\ldots,q-1$.

Let $Y=\sum _{s\in Z/p}\sum_{a=0}^{p-1}
    \tau^Y_{a,s}1_{s}E^{(a)}F^{(a)}$ be in $ Z(A)$. Then for any $s\in Z/p$,
$$
1_{s}E^{(1)} Y=Y 1_{s}E^{(1)}
$$
From here by direct computation one can see that
 for any $0\leq a\leq p-2$,
$$
[a-s]\tau^Y_{a+1,s+2}=[a+1](\tau^Y_{a,s}-\tau^Y_{a,s+2}).
$$
This implies that if $Y\in \Omega_{b}$ then
$\tau^Y_{ b,s}$ doesn't depend on $s$ and we denote it with
$\tau^Y_{b}$.
In particular, $X^{b}$ is of this type, moreover $\tau^X_{b}\neq 0$
and therefore, there exists $r\in k$ such
that  if $b>0$ then
$Y-rX^{b-1}\in \Omega_{b-1}$ and if $b=1$ then $Y=r X^0=r \1$.
The proposition follows by induction.\end{proof}

\begin{ssec}\label{pn} {\bf Proof of lemma \ref{Jact}}\qua
    Now we continue with the proof of the lemma \ref{Jact}.
Observe that since $[N_i^-+N_i^+]=0$,
$\hat{J}([N_i^-])=-\hat{J}([N_i^+])$, so we will
compute $\hat{J}([N_i^+])$.
From the expressions in \ref{eval} one obtains:
\beq
&&P_j=1_{-2j}+1_{2j+2}+\sum_{s\in Z/p}\sum_{a=1}^{p-1}\tau
^j_{s,a}1_{s}E^{(a)}F^{(a)},\,0\leq j\leq q-1\\
&&P_q=1_{1}+\sum_{s\in Z/p}\sum_{a=1}^{p-1}\tau
^q_{s,a}1_{s}E^{(a)}F^{(a)},\\
&&N_j=\sum_{a=0}^{p-1}\nu ^j
_{-2s,a}1_{-2s}E^{(a)}F^{(a)},\;
N_j^+=\sum_{s=j+1}^{p-1-j}\sum_{a=0}^{p-1}\nu ^j
_{-2s,a}1_{-2s}E^{(a)}F^{(a)},
\eeq
where $\nu ^j_{-2s,0}=0$ and $\nu ^j_{-2s,p-1}=
    (v-v^{-1})[2j+1]^2$.
Given $i,a$  such that
$0\leq i\leq q-1$, $0\leq a\leq p-1$ and given $s\in Z/p$,
let $\bar{\nu} ^i_{-2s,a}\in k$ are the coefficients of the expansion
of $J(N_i^+)$ in terms of the basis $1_{-2s}E^{(a)}F^{(a)}$, i.e.\
$$
J(N_i^+)=S\circ J(N_i^+)=\sum _{n,m}\lambda(\beta _n
N_i^+\alpha _m)S(\alpha _n\beta _m)=
\sum _{s\in Z/p}\sum_{a=0}^{p-1}\bar{\nu} ^i_{-2s,a}1_{-2s}E^{(a)}F^{(a)}.
$$
Substituting here the expression for the $R$--matrix and for $N_i^+$ 
we obtain
$$
\bar{\nu} ^i_{-2s,a}=v^{a(a+1)+2as}\{a\}^2\left[\begin{array}{c}
p-1\\a\end{array}\right]^2
\sum_{l=i+1}^{p-1-i}v^{2l(a-2s-1)}\nu^i_{-2l,p-1-a}.
$$
In particular
$\bar{\nu} ^i_{-2s,p-1}=0$ and
$$
\bar{\nu} ^i_{-2s,0}=-(v-v^{-1})[2i+1]^2\frac{[(2s+1)(2i+1)]}{[2s+1]}.
$$
Then the lemma follows from  proposition \ref{poly}
and the expression for $P_s$.
\end{ssec}
\begin{ssec}\label{dotn}
For any $0\leq i,j\leq q-1$, let
$
\omega_{i,j}=\frac{[(2j+1)(2i+1)]}{[2j+1]}$. Let also
$$
\dot{N}_i=\frac{N_i}{(v-v^{-1})[2i+1]^2}\;\mbox{ 
and }\;
\dot{N}_i^\pm=\frac{N_i^\pm}{(v-v^{-1})[2i+1]^2}.
$$
Observe that $[\dot{N}_0^-]=[\Lambda]$. 
Since $\hat{J}$ is injective, the $(q-1)\times (q-1)$ matrix $\omega$
is nondegenerate and the proposition above implies that
$$
\hat{J}([\dot{N}_j^{-}])=\sum_{i=1}^{q-1}\omega _{ji}[P_i]\;\mbox{ 
and }\;
\hat{J}([P_j])=\gamma_p\sum_{i=1}^{q-1}(\omega^{-1}) _{ji}[\dot{N}^-_i].
$$
\end{ssec}
 \begin{prop}\label{kg} \begin{itemize}
\item[\rm(a)] $ \sigma_{ij}^k=\sigma
([\dot{N}_i^{-}],[\dot{N}_j^{-}],[P_k])=\lambda (\dot{N}_k^-)
\sum_{s=0}^{q-1}\omega_{is}\omega_{js}\omega^{-1}_{sk}$,
and  $\sigma(a,b,c)=0$ for any other triple of generators $a,\,b,\, c$;
\item[\rm(b)]     $\sigma_{ij}^k/\lambda (\dot{N}_k^-)=1$ if
    all of the following four conditions are satisfied:
    $$
    i+j+k\leq p-2,\quad i+j-k\geq 0,\quad k+i-j\geq 0,\quad k+j-i\geq 0.
    $$
Otherwise $\sigma_{ij}^k=0$.
\end{itemize}
\end{prop}

\proof[Proof of (a)]
Lemma \ref{J} allows as to express the $\star $ product in the
following way:
\beq
&&[\dot{N}_{i}^-]\star [\dot{N}_{j}^-]=
  \gamma_p^{-1} \hat{J}(\hat{J}([\dot{N}_{i}^-])\hat{J}([\dot{N}_{j}^-]))=
\sum_{k,s=0}^{q-1}\omega_{is}\omega_{js}(\omega^{-1})_{sk}[\dot{N}_{k}^-];\\
&&[P_{i}]\star [P_{j}]= \gamma_p^{-1}
\hat{J}(\hat{J}([P_{i}])\hat{J}([P_{j}]))=0;
\eeq
This implies that $\sigma(a,b,c)=0$ if all three elements are
of the type $\dot{N}_i$, or if only one of them is such. For the only
nonzero case
we obtain
$$
\sigma
([P_k],[\dot{N}_i^{-}],[\dot{N}_j^{-}])=\lambda(P_k,(\dot{N}_i^{-}\star
\dot{N}_j^{-}))=
\lambda (\dot{N}_k^-)
\sum_{s=0}^{q-1}\omega_{is}\omega_{js}(\omega^{-1})_{sk}.\eqno{\qed}
$$

\begin{proof}[Proof of (b)]
 Using that for any primitive  $p$-th root of
unity $v$ and any $a\in Z/p$,
$$
\sum_{s=0}^{q-1}v^{a(2s+1)}=\frac{p\,\delta_{a,0}-v^{-a}}{1+v^{-a}},
$$
one obtains that
$$
\sum_{i=0}^{q-1}[(2j+1)(2i+1)][(2i+1)(2k+1)]=-\frac{p}{(v-v^{-1})^2}\delta_{
j,k}.
$$
Hence
$$
(\omega^{-1})_{i,j}=-\frac{(v-v^{-1})^2}{p}[2i+1][(2i+1)(2j+1)],
$$
and
$$
\frac{\sigma_{ij}^k}{\lambda
(\dot{N}_k^-)}=-\frac{(v-v^{-1})^2}{p}
\sum_{s=0}^{q-1}\frac{[(2i+1)(2s+1)][(2j+1)(2s+1)][(2k+1)(2s+1)]}{[2s+1]}.
$$
Substituting above the expression
$$
\frac{[(2i+1)(2s+1)]}{[2s+1]}=\sum_{l=0}^{2i}v^{2(i-l)(2s+1)},
$$
and expanding we obtain that
\beq
&&p\frac{\sigma_{ij}^k}{\lambda
(\dot{N}_k^-)}\\
&&=\sum_{l=k-i-j}^{k+i-j}\sum_{s=0}^{q-1}(v^{2l(2s+1)}+v^{-2l(2s+1)})-
\sum_{l=j+k-i+1}^{j+k+i+1}\sum_{s=0}^{q-1}(v^{2l(2s+1)}+v^{-2l(2s+1)})\\
&&=p\left(\sum_{l=k-i-j}^{k+i-j}\delta_{\bar{l},0}-
\sum_{l=j+k-i+1}^{j+k+i+1}\delta_{\bar{l},0}\right),
\eeq
where $\bar{l}=\mbox{Mod}(l,p)$. This completes the proof of the
proposition.\end{proof}
Observe that the proof of proposition \ref{kg} above 
 imply:
 \begin{corr} The subalgebra of $\hat{Z}_{\star}(A)$
spanned by  $[\dot{N}_j^-]$, $0\leq j\leq (q-1)$ is 
isomorphic to the fusion algebra ${\cal F}_p$ of
the semisimple quotient of the representation category of $A$
defined in \ref{fp}. 
\end{corr}

Finally we  can
describe all elements in $\T_{Z}$.
\begin{theo} \label{sl2} $\T_{Z}$ consists of the multiples of 
$[1],\,[\Lambda] ,\, \sum_{j=0}^{q-1}[2j+1][\dot{N}_j^-]$ and $[P_0]$.
Moreover, $\hat{J}$ sends bijectively $\T_{Z}$ into itself.
\end{theo}
\begin{proof}
 Suppose that
$[z]=\sum_{i=0}^{q-1}x_i[P_{i}]+\sum_{i=0}^{q-1}y_i[\dot{N}_{i}^-]
$.
 According to \ref{prop2}  $[z]\in\T_Z$ if and only if
 for any $[a],[b],[c]\in \hat{Z}(A)$,
$
\sigma (zc,za,b)=\sigma (zc,a,zb)
$.
Replacing here all possible choices of $a,b,c$ we obtain that this
condition is equivalent to the following system of equations for the
coefficients
$x_i,y_i$:
\beq
&\text{(i)}&y_iy_k\sigma_{ik}^j=y_jy_k\sigma_{jk}^i;\\
&\text{(ii)}&y_k(x_i-x_j)\sigma_{jk}^i=y_ix_k\sigma_{ji}^k;\\
&\text{(iii)}&x_k(x_i-x_j)\sigma_{ji}^k=0;\\
&\text{(iv)}&y_ix_k\sigma_{ik}^j=y_jx_k\sigma_{jk}^i;\\
&\text{(v)}&x_k(x_i-x_j)\sigma_{jk}^i=0,
\eeq
for any $0\leq i,j,k\leq q-1$.
Now we want to show that $[z]$
is contained either in the span of the $[P_i]$'s or in the span of the
$[\dot{N}_i^-]$'s. Observe that
$\sigma_{0,j}^k=\delta_{kj}\lambda(\dot{N}_j^{-})$.
Hence equations (v) and (ii) with $j=0$ become
$$
x_i(x_i-x_0)=0\qquad y_ix_0=0.
$$
Therefore either $x_i=0$ for any $i$ or  $y_i=0$ for any
$i$. Suppose now that we are in the case when $y_i=0$ for any
$i$ and let ${\cal I}\neq \emptyset$ be  the subset of
indices such that $x_i\neq 0$. Then,  condition (v) implies that
\begin{itemize}
\item[(a)] $0\in {\cal I}$;
\item[(b)]  for any other $i\in {\cal I}$, $x_0=x_i$;
\item[(c)]   if $\sigma_{jk}^i\neq 0$ and
two of the indices $i,j,k$ are in ${\cal I}$, then the third one must be in 
${\cal I}$ as well.
\end{itemize}
 Moreover, any subset ${\cal I}$ which satisfies these conditions corresponds to
 a solution of the form $[z_{\cal I}]=\sum_{i\in {\cal I}}[P_i]$.
 In particular,  since $\sigma_{0,0}^k=\sigma_{0,k}^0=\delta_{k,0}$,
 ${\cal I}=\{0\}$ ($[z_{\cal I}]=[P_0]$) gives  a solution of the problem.

Suppose now that $i\in {\cal I}$ and $i\neq 0$. Since
 $\sigma_{ii}^1\neq 0$
(\ref{kg} (b)) it follows that  $1$ should be in ${\cal I}$ as well. 
But if $1,j\in {\cal I}$ where $j\leq q-2$, then  $j+1\in {\cal I}$ (since
$\sigma_{j,1}^{j+1}\neq 0$).
Hence, if ${\cal I}$ contains one nonzero index, it must contain all 
indices, i.e.\ 
 ${\cal I}=\{0,1, \ldots, q-1\}$ and $[z_{\cal I}]=[\1]$.

Suppose now that $x_0=0$ and  ${\cal I}\neq \emptyset$ is  the subset of
indices such that $y_i\neq 0$  i.e.\  $[z]=\sum_{i\in {\cal I}}y_i[\dot{N}_i^-]$.
From  \ref{kg} (b) it follows that 
$\sigma_{jk}^i=\lambda(\dot{N}_i)\epsilon_{ijk}$ where 
$\epsilon_{ijk}$ is symmetric with respect to the three indices. Then 
 equation (i)  becomes:
 $$y_k(\lambda(\dot{N}_j)y_i-\lambda(\dot{N}_i)y_j)\epsilon_{ijk}=0.
 $$
 In particular for  $i=0$ and $j=k$ we have 
 $y_k(\lambda(\dot{N}_k)y_0-y_k)=0$. 
Hence ${\cal I}$ satisfies the conditions (a)--(c) above and 
therefore either ${\cal I}=\{0\}$ or  ${\cal I}=\{0,1, \ldots, q-1\}$ and 
$y_k=\lambda(\dot{N}_k)y_0$  for any 
$k\neq 0$. 
 The   corresponding solutions for $[z]$ are
 $[z]=y_0[N_0^-]=y_0\gamma_p [J(\1)]$ and
 $$
 [z]=y_0\sum_{j=0}^{q-1}[2j+1][\dot{N}_j^-]=-\frac{y_0
\hat{J}([P_0])}{(v-v^{-1})^2p^2}. 
$$
This completes the proof of the
theorem.\end{proof}

\begin{ssec} {\bf The $sl(2)$ HKR--type invariants}
       
    We remind that $\T_s$ denotes 
    the subset of  elements in $\hat{Z}^S(A)$ which 
    define invariants of links under the band-connected sum of two 
    distinct components. 
    Then $\T^3\subset \T\subset\T_s$. 
  We can not offer a
   way to calculate the elements in $\T$ and even less a way to study
   its maximality, i.e.\ if it coincides 
    with $\T_s$. But since $\T_Z\supset\T$, a
    hypothetical search for the elements in $\T$ could 
    start by calculating the elements in  $\T_Z$ as it has been 
    done above for the  $sl(2)$ case. The surprise is that $\T_Z$ is 
    already very restrictive:  up to multiplication by an element 
    in $k$,
      it consists of   four elements and, 
    using proposition \ref{rts} in the appendix,
    we  see that three of them are  in $\T_s$:
    \beq
    &&[z_H]=[\1] \;\mbox{ gives the Hennings invariant};\\
     &&[z_{RT}^*]=[P_0] \\
     &&[z_H^*]=[\Lambda ] \;\mbox{ gives the trivial invariant (equal to 1 for
     any manifold)};\\
      & &[z_{RT}]=-\frac{[J(P_0)]}{(v-v^{-1})^2p^2}=
             \sum_{j=0}^{q-1}[2j+1][\dot{N}_j^-]
     \;\mbox{ gives the RT--invariant};
     \eeq
     So, it seems reasonable to make the following conjecture:

     \begin{conj}\label{ttz} If $A$ is a finite dimensional, unimodular, ribbon,
	 $\Lambda$--factor\-iz\-able algebra, then $\T_Z=\T$.
	 \end{conj}
	If the conjecture holds then
	 $sl(2)$
     produces exactly four HKR--type invariants, all normalizable
     to 3--manifold invariants. Moreover, since
    $
     [P_0 z_{RT}]=-p^2 [\Lambda ]$ and $[P_0\star z_{RT}]=[1]
     $, proposition \ref{3factfact}
     implies:
 \begin{corr}\label{hrt}$
     \Z_{[ z_H]}^{\partial}(\partial M)=
     \Z_{[ z_{RT}]}^{\partial}(\partial
M)\Z_{[z_{RT}^*]}^{\partial}(\partial M).
$\end{corr}

To support the conjecture, we  show 
that the statement of corollary \ref{hrt} holds for the
values of the three invariants for $S^2\times S^1$ and
the Lens spaces. Directly from the definition for $S^2\times S^1$ we 
have:
\beq
&&\Z_{[z_H]}^{\partial}( S^2\times S^1)=\lambda(1)=0;\;\\
&&\Z_{[z_{RT}^*]}^{\partial}( S^2\times
S^1)=\lambda(P_0)=0;\;\\
&&\Z_{[z_{RT}]}^{\partial}( S^2\times
S^1)=\lambda(z_{RT})=\sum_{j=1}^{q-1}=[2j+1]^2.
\eeq
Observe that 
$\Z_{[z]}^{\partial}(L(1,n))=\lambda(z\theta^n)/\lambda(z\theta)$. 
Then from \ref{th} one obtains  that
$
[\theta]^n=\sum_{j=0}^{q-1}v^{2nj(j+1)}([P_j]-n\,p(v-v^{-1})[2j+1][\dot{N}_j^-])
$.
Hence,
\beq
&&\lambda(z_{RT}^*\theta^n)=
-pn(v-v^{-1});\;\\
&&\lambda(z_{RT}\theta^n)=
\sum_{j=0}^{q-1}v^{2nj(j+1)}[2j+1]^2;\\
&&\lambda(\theta^n)=
-pn(v-v^{-1})\sum_{j=0}^{q-1}v^{2nj(j+1)}[2j+1]^2.\\
\eeq
Therefore the statement of corollary \ref{hrt} holds for the
values of the three invariants for 
the Lens spaces as well.
\end{ssec}

\section{Questions}

\begin{ssec}  If the conjecture \ref{ttz}
     is false, this would  imply that the
 condition $[z]\in\T$ in theorem \ref{theo1} is too strong and needs to be weakened.
Then one may ask if it can be replaced with $[z]\in \T_Z$.
\end{ssec}

 \begin{ssec} In the case of the quantum $sl(2)$ we saw that the
 fusion algebra of the semisimple quotient of the representation
 category is a subalgebra of $\hat{Z}^S_{\star}(A)$ generated by nilpotent
 elements.
 What is in general the relationship between $
\hat{Z}^S_{\star}(A)$ and the representation theory of $A$?
 \end{ssec}

 \begin{ssec} Observe that if   the Hopf algebra is
triangular, then $\T^3=\{ X[\Lambda]\mid X\in k\}$, i.e.\ such algebra
doesn't produce nontrivial 3--manifold invariants. On another hand if
$\T^3=\T^4$ (i.e.\ any 4--invariant is normalizable to a 3--manifold
invariant) then $\T^2= \{ X[\Lambda]\mid X\in k\}$, i.e.\ such algebra
doesn't produce nontrivial  invariants of 2--complexes. This seems to be
the example of the quantum $sl(2)$. It would be interesting to know
if there exists a Hopf algebra for which $\T^4$ doesn't reduce to 
$\T^2$ or to $\T^3$ and if a similar algebra exists, one may ask if 
as $\star$--monoid $\T$ is generated by $\T^2$ and $\T^3$. This is 
related to the following purely topological question:
\end{ssec}

\begin{ssec} Let $(M,P)$ and $(M',P')$ be two 4--thickenings such that
    ${\rm index}(M)={\rm index}(M')$,
    $P$ is 2--equivalent to $P'$ and $\partial M$ is diffeomorphic to 
    $\partial M'$. Then is it true that $M$ is diffeomorphic to $M'$? Is
    $M$ 2--equivalent to $M'$? The results in \cite{Q:dual} seem to support 
    the affirmative answer.
 \end{ssec}

\section{Appendix: The Reshetikhin--Turaev $sl(2)$--invar-iant  as HKR--type
invariant}

Before starting working on this project, the first author  asked
T.Kerler why the Reshetikhin--Turaev $sl(2)$--invariant
is  a HKR--type invariant.
For  completeness we give here
Kerler's explanation and the evaluation of the corresponding
trace element $[z_{RT}]\in \hat{Z}(A)$. In somewhat different form
this evaluation has been done in \cite{K:g}.

We use the definition of the Reshetikhin--Turaev invariant as
given in \cite{RT}. But since the
definition of the quantum $sl(2)$ here  is slightly different
from the one in \cite{RT}, the reader is
referred to the work of Gelfand and Kazhdan \cite{GK}
for the proof that, the full linear category generated from
the ``small'' representations used below,
satisfies the requirements in paragraph 3.1 of \cite{RT}.

Let $A$, $k$ and $g$ be as in \ref{defn}. For any finite dimension
left $A$--module $V$ define the dual representation $V^*$ of $V$
to be representation with linear space $Hom(V,k)$ and action of
$a\in A$ given by $S(a)^*$.
Define also the quantum trace $\tr_V\co A\rightarrow k$ of $V$ to be
$$
\tr_V(a)=\sum _{i=1}^{dim(V)}e_i^*(gae_i)\;\mbox{for any $a\in A$},
$$
where $\{e_i\}_{i=1}^{dim(V)}$ is a basis for $V$ and $\{e_i^*\}_{i=1}^{dim(V)}$
its the dual basis for $V^*$.

\begin{prop}\label{zv} For any finite dimensional
left $A$--module $V$ there exists $z_V\in Z(A)$ such that for any
$a\in A$,
$\tr_V(a)=\lambda(g^2z_Va)$.
\end{prop}
\begin{proof}
First observe that for any $a,b\in A$,
$$
\tr_V(ab)=\sum _{i=1}^{dim(V)}e_i^*(gabe_i)=\sum _{i=1}^{dim(V)}e_i^*(bgae_i)=
\tr_V(bS^2(a)).
$$
Now, since $g$ is invertible, from \ref{print} (d) it follows that there
exists an
element $z_V\in A$
such that $\tr_V(a)=\lambda(g^2z_Va)$. Moreover,  \ref{unimodular} and
\ref{ribbon} imply that for any $a,b\in A$,
\beq
&&\lambda(g^2(z_Va-az_V)b)=\lambda(g^2z_Vab-S^4(a)g^2z_Vb)\\
&&=\lambda(g^2z_Vab-g^2z_VbS^2(a))=\tr_V(ab-bS^2(a))=0.
\eeq
Then the statement
follows from \ref{print} (d).\end{proof}

\begin{ssec} \label{repr}Let $\Sigma=\{0,1, \ldots,q-1\}$ and
let $V_n$, $n\in\Sigma$ be the simple left $A$-module
with highest weight $2n$. Then $V_n$ has a basis
$\{e^n_{i}\}_{i=-n}^n$ and the action of the algebra generators is as 
follows:
\beq
&&1_{2i-2} F^{(1)}e^n_{i}=\left\{\begin{array}{lc} 0&\mbox{if}\; i=-n\\
                          
                          e^n_{(i-1)}&\mbox{otherwise}\end{array}\right. 
                          \\
&& 1_{2i+2} E^{(1)}e^n_{i}=\left\{\begin{array}{lc} 0&\mbox{if}\; i=n\\
        \mbox{$[n+i+1]$} e^n_{(i+1)}&\mbox{otherwise}\end{array}\right. 
        \\
&&    1_ce^n_{i}=\left\{\begin{array}{lc}  e^n_{i}&\mbox{if}\;
    c=2i\;(mod \;p)\\
     0&\mbox{otherwise}\end{array}\right.
\eeq
When it is clear which one is the representation, we will use $e_{i}$
instead of $e^n_{i}$. Moreover, $d_n=2n+1$ will denote the dimension
of $V_n$ and $z_n=z_{ V_n}$.

Given a sequence ${\bf i}=(i_1,i_2, \ldots,i_k)$ of elements in
$\Sigma$,  define 
$$V({\bf i})=V_{i_1}\otimes V_{i_2}\otimes \ldots 
\otimes V_{i_k}\;\mbox{ and }\;
r({\bf i})=tr_{V({\bf i})}(id).
$$
Observe that $r(n)=[2n+1]$ and since $g$ is a group-like element, $r({\bf
i})=\prod_{s=1}^kr(i_s)$.
\end{ssec}
\begin{ssec}\label{fp} As it is shown in \cite{GK}, the full linear category ${\cal C}_p$ generated by $V_n$, 
$n\in\Sigma$, is equivalent to the semisimple quotient of the category 
of integral representation of $A$, and this equivalence 
induces a braided monoidal structure on ${\cal C}_p$. In particular 
there is a product structure on ${\cal C}_p$ given by
$$
V_i\diamond V_j=\oplus_{s\in \Sigma }k^{\epsilon_{ij}^s}\otimes V_s.
$$
The essence of this product structure is encoded in the fusion 
algebra ${\cal F}_p$ which is defined as the vector space $Z[x_0, 
x_1, \ldots,x_{q-1}]$ 
and product structure given by 
$$ x_ i\diamond  x_j=\sum_{s\in \Sigma }\epsilon_{ij}^s\,x_s,
$$
for any $ i, j\in \Sigma$.
The (non negative) integers $\epsilon_{ij}^s$ are called the fusion 
coefficients of ${\cal C}_p$. The fusion coefficients for the quantum 
$sl(2)$ have been calculated in \cite{RT,GK} and are the following:
$\epsilon_{ij}^s=1$ if
    all of the following four conditions are satisfied
    $$
    i+j+s\leq p-2,\quad i+j-s\geq 0,\quad s+i-j\geq 0,\quad s+j-i\geq 0,
    $$
and $\epsilon_{ij}^s=0$ otherwise. 
\end{ssec}
\begin{ssec}
Given an oriented $k-l$ tangle $T$,
represented with a tangle diagram, one associates to the incoming
and the outgoing ends of $T$ the sequences
$\underline{\epsilon}=\{\epsilon_1, \ldots
,\epsilon_k\}$ and $\overline{\epsilon}=\{\epsilon^1, \ldots
,\epsilon^l\}$ where $\epsilon_i=1$ ($\epsilon^i=1$) if in a neighborhood
of the point
the tangle component points down and $\epsilon_i=-1$ ($\epsilon^i=-1$)
otherwise.

A {\em coloring} ${\bf n}=(n_1,n_2,\ldots,n_m)\in\Sigma^{\times m}$ 
of an oriented $k-l$ tangle $T$ with $m$
components, is a
map which associates to the $i$-th connected component of $T$ an element 
$n_i\in \Sigma$. A coloring of the tangle induces  colorings $\underline{ i}({\bf n})=
\{i_1, i_2, \ldots,i_k\}$ and $\overline{ i}({\bf n})=\{i^1, i^2, \ldots,i^l\}$
of the incoming and the outgoing ends of the tangle.

The colored tangles form a category $\h$ with objects the set
${\cal S}$ of sequences $\{(\epsilon_s, i_s)\}_{s=1}^k$, where
$\epsilon _s=\pm 1$ and $i_s\in \Sigma$.
If $\eta,\eta'\in {\cal S}$ then a morphism $\eta\rightarrow
\eta'$ is a colored tangle considered up to isotopy such that the
sequence of signs and colors of the outgoing ends is equal to $\eta$
and the one of the incoming ends is equal to $\eta'$ (This is not a
mistake. While in the HKR framework we were multiplying the algebra
elements on the right,
in the Reshetikhin--Turaev framework one considers the left action of the
algebra on a
representation and this leads to the necessity of reversing the idea of
incoming and outgoing). The composition of two tangles $T'\circ T$ is
obtained by placing  $T'$ on the top of $T$ and gluing the ends.
The category can also be
provided with tensor product by defining $T'\otimes T$ to be the
tangle obtained by placing $T'$ to the left of $T$.
\end{ssec}

\begin{ssec}\label{defrt}
Theorem 2.5 in \cite{RT} states that there exists
a unique covariant functor $F\co \h\rightarrow Rep\, A$  such that
for any object $\eta$ in $\h$,
$
F(\eta)=V^{\epsilon_1}_{i_1}\otimes V^{\epsilon_2}_{i_2}\otimes \ldots
V^{\epsilon_k}_{i_k}$,
where $V_n^1=V_n$ and $V_n^{-1}=V_n^*$.
Moreover, $F$ preserves the tensor product and if $F(T;{\bf n})$ denotes
the value of $F$ on an oriented tangle $T$ with coloring ${\bf n}$, on the
elementary colored tangles  presented in figure \ref{fig.tangelem}
this value is as follows:
\beq
&&F(b1;i)=id_{V_i},\quad F(b2;i)=id_{V_i^*}\\
&&F(d1;i,j)\co x\otimes y\rightarrow \sum_{n}\beta_n.y\otimes
\alpha_n.x\co V_i\otimes V_j\rightarrow V_j\otimes V_i;\\
&&F(d2;i,j)\co x\otimes y\rightarrow \sum_{n}S(\alpha_n).y\otimes
\beta_n.x\co V_i\otimes V_j\rightarrow V_j\otimes V_i,\\
&&F(e1;i)\co x\otimes y\rightarrow x(y)\co V_i^*\otimes V_i\rightarrow k;\\
&&F(e2;i)\co y\otimes x\rightarrow x(g.y)\co V_i\otimes V_i^*\rightarrow k;\\
&&F(f1;i)\co 1\rightarrow \sum_{k=1}^{d_i}e_k\otimes
e_k^*\co k\rightarrow V_i\otimes V_i^*;\\
&&F(f2;i)\co 1\rightarrow \sum_{k=1}^{d_i}e_k^*\otimes
g^{-1}.e_k\co k\rightarrow V_i^*\otimes V_i,
\eeq
where with ``.'' denotes the left action of $A$ on the
corresponding left $A$--module.

Let $L$ be a link with $m$ components. Fix an orientation of $L$ and
define $
\{L\}=\sum_{{\bf n}}r({\bf n})F(L;{\bf n})
$, where the sum is over all possible colorings ${\bf n}=\{n_1,\ldots,n_m\}$
of $L$.
Then theorem 3.3.2 in \cite{RT} states that
$
\{L\}$
doesn't depend on the orientation of the components of $L$. Moreover,
$\{L\}$ is an invariant of the link under isotopy and under taking
the band connected sum of two different components.
\end{ssec}

\begin{prop} \label{rts} $\{L\}=\Z_{[z_{RT}]}(L)$, where
  $z_{RT}=\sum_{n=0}^{q-1} r(n) z_n$.  In particular, $[z_{RT}]\in \T_s$.
    \end{prop}
\begin{proof}
 We can represent $L$ as the closure of a braid $B$
 on $k$ strings oriented downwards as in the example in figure
 \ref{braid}.
       \begin{figure}[ht!]
\setlength{\unitlength}{1cm}
\begin{center}
\begin{picture}(5,4)
\put(0,0){\fig{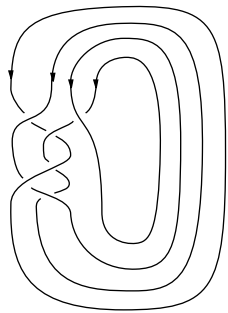}}
\end{picture}
\end{center}
\caption{Presenting a link as the closure of a braid}
\label{braid}
\end{figure}
 Let $\sigma$ be the underlying permutation of $B$, i.e.\ the boundary
 of the $i$-th
 component of $B$ consists of the $i$-th incoming and the
 $\sigma(i)$-th outgoing ends (counted from the left to the right).
 Then $\sigma$ is the product of $m$ cycles:
 $$\sigma=(j^1_1,j^1_2, \ldots j^1_{s_1})\ldots(j^m_1,j^m_2, \ldots j^m_{s_m}).
 $$
 For the example of figure \ref{braid}, $\sigma=(1)(2,3,4)$.
 Let
 $\Z(B)=\sum_{i}a_{1,i}\otimes a_{2,i}\ldots\otimes a_{k,i}$
 be the element in  $A^{\otimes k}$ as defined in \ref{defmap}.
  Let also
 $$\sum_{j}c_{1,j}\otimes c_{2,j}\ldots\otimes c_{m,j}=
 \sum_{i}(a_{j^1_1,i}ga_{j^1_2,i}\ldots ga_{j^1_{s_1},i})\otimes \ldots
 \otimes (a_{j^m_1,i}ga_{j^m_2,i}\ldots ga_{j^m_{s_m},i}).
 $$
Then from \ref{defrt} it follows that 
for any coloring ${\bf n}=\{n_1, n_2, \ldots,n_m\}$ of $L$,
 $$
 F(L;{\bf n})= \sum_{j}\tr_{V_{n_1}}(c_{1,j}) \ldots
 \tr_{V_{n_m}}(c_{m,j})=\sum_{j}\lambda(g z_{n_1}c_{1,j}g) \ldots
 \lambda(gz_{n_m}c_{m,j}g).
 $$
 Here we have used   the fact
 that for any $a,b\in A$ and  $-n\leq s,l\leq n$,
  $$\sum_{i=-n}^{n}e_{l}^*(a.e_{i})e_{i}^*(b.e_{s})=e_{l}^*(ab.e_{s}).
  $$
  Making the confrontation with the expression for $\Z$ in 
  \ref{defmap},
  we see that
  $F(L;{\bf n})$ $=\Z(L)(z_{n_1}, \ldots,z_{n_m})$. The
 statement of the proposition follows by linearity.\end{proof}

 \begin{prop} For any $0\leq n\leq q-1$, $[z_n]=[\dot{N}_n^-]$.
     \end{prop}
\begin{proof}
     From \ref{th} it follows that
     $z_n=\sum_{i=0}^{q-1}(x_iP_i+y_i\dot{N_i}^-+w_i\dot{N_i})+x_qP_q$.
     Let $0\leq j\leq q$ and
     $a_j=1_{-2j}E^{(p-1)}F^{(p-1)}$. Then  \ref{repr} implies that
     $\tr_{V_n}(a_j)=0$ for any $j$. On
     another hand, from the expressions for $P_j$ and $N_j$ in \ref{pn}
     it follows that
     $$\lambda(gz_n a_jg)=v^{2j}x_j.
     $$
     Hence $x_j=0$ for any $0\leq j\leq q$. On another hand, for every $0\leq
     j\leq n$,
     $\tr_{V_n}(1_{-2j})= v^{2j}$ and from \ref{pn}  it follows that
     $$
     \lambda(gz_n
     1_{-2j}g)=v^{4j}\sum_{i=0}^{q-1}(y_i\lambda(\dot{N_i}^-1_{-2j})+
     w_i\lambda(\dot{N_i}1_{-2j})).
 $$
Hence we obtain the following system of $q+1$ equations for the coefficients
$y_i,w_i$:
\beq
&&\sum_{i=j}^{q-1}y_i+\sum_{s=0}^{q-1}w_s=1,\;0\leq j\leq n;\\
&&\sum_{i=j}^{q-1}y_i+\sum_{s=0}^{q-1}w_s=0,\;n+1\leq j\leq q;
\eeq
The solution is $y_i=\delta_{i,n}$ and $\sum_{s=0}^{q-1}w_s=0$. 
Hence $z_n=\dot{N}_n^-+\sum_{s=0}^{q-1}w_s\dot{N}_s$ and 
$[z_n]=[\dot{N}_n^-]$.
\end{proof}
As a consequence of the last two propositions it follows that
$$[z_{RT}]=\sum_{n=0}^{q-1}[2n+1][\dot{N_i}^-].$$

\Addresses
\end{document}